\documentclass[a4paper, reqno]{amsart}
\usepackage{amsmath, amssymb, amscd}

\numberwithin{equation}{section}
\theoremstyle{plain}
\newtheorem{thm}{Theorem}[section]
\newtheorem{cor}[thm]{Corollary}

\newtheorem{lem}[thm]{Lemma}
\newtheorem{prop}[thm]{Proposition}

\theoremstyle{remark}
\newtheorem{rem}[thm]{Remark}

\title{An expansion formula for the Askey-Wilson function}
\author{Jasper V. Stokman}
\address{Jasper V. Stokman,
KdV Institute for Mathematics, Universiteit van Amsterdam,
Plantage Muidergracht 24, 1018 TV Amsterdam, The Netherlands.}
\email{jstokman@science.uva.nl}
\date{May 11th, 2001\\
\indent 2000 {\it Mathematics Subject Classification}. Primary 33D45;
Secondary 33D80.}

\begin{document}

\keywords{}

%%%%%%%%%%%%%%%%%%%%%%%%%%%%%%%%%%%%%%%%%%%%%%%%%%%%%%%%%%
%%                                                      %%
%%             Abstract                                 %%
%%                                                      %%
%%%%%%%%%%%%%%%%%%%%%%%%%%%%%%%%%%%%%%%%%%%%%%%%%%%%%%%%%%
\begin{abstract}
The Askey-Wilson function transform is a $q$-analogue of the Jacobi
function transform with kernel given by an explicit 
non-polynomial eigenfunction of
the Askey-Wilson second order $q$-difference operator.
The kernel is called the Askey-Wilson function.
In this paper an explicit expansion formula for 
the Askey-Wilson function in terms of Askey-Wilson polynomials
is proven. With this expansion formula at hand,
the image under the Askey-Wilson function transform of
an Askey-Wilson polynomial multiplied by an analogue
of the Gaussian is computed explicitly.
As a special case of these formulas a $q$-analogue
(in one variable) of the Macdonald-Mehta integral is obtained,
for which also two alternative, direct proofs are presented.
\end{abstract}
\maketitle
\section{Introduction}
The Macdonald polynomials and their orthogonality relations
have an harmonic analytic interpretation on
quantum compact Riemannian symmetric spaces, see e.g. Noumi \cite{N}.
In particular, the spherical
Fourier transform on the quantum $SU(2)$ group
can be identified with the 
polynomial Askey-Wilson transform,
which is the transform naturally associated to the orthogonality relations
of the Askey-Wilson polynomials (see Koornwinder \cite{K}).

In the non-compact set-up,
only harmonic analysis on the quantum $SU(1,1)$ group
has been well understood by now, see e.g. \cite{KMU} and \cite{KS1}.
This has led to the study of an explicit
generalized Fourier transform in \cite{KS1} and \cite{KS2},
called the Askey-Wilson function transform. The kernel of this
transform is called the Askey-Wilson function. It is a non-polynomial
eigenfunction of the Askey-Wilson second-order
$q$-difference operator, given explicitly as a
very-well-poised ${}_8\phi_7$ series.

On the other hand, Cherednik \cite{C} discussed several types of
difference Fourier transforms, which are naturally related to 
the spectral theory of Macdonald polynomials.
Cherednik \cite{C} showed that a particular theta-function
plays a role in the theory of the difference transforms which is
similar to the role of the Gaussian in the theory of Hankel
transforms, see also \cite{CM}. 
This led to explicit formulas for
the image under the difference Fourier transforms
of a Macdonald polynomial multiplied by (the inverse of) the
analogue of the Gaussian.
Furthermore, in the difference set-up Cherednik \cite{C} defined certain
``non-polynomial'' spherical functions as explicit
series expansions in terms of ``polynomial'' spherical functions
(=Macdonald polynomials), which seems to be
purely a quantum phenomenon.

The purpose of the present paper is to incorporate
the above mentioned ideas and constructions of Cherednik into 
the theory of the polynomial Askey-Wilson transform, and into the
theory of the Askey-Wilson
function transform. We first give an explicit expansion formula
for the Askey-Wilson function in terms of
Askey-Wilson polynomials. This expansion formula 
provides an explicit link between Cherednik's \cite{C} construction of
non-polynomial eigenfunctions of $q$-difference operators
with the
constructions of Suslov \cite{S89}, Ismail \& Rahman \cite{IR}
using the theory of basic hypergeometric series.
We introduce the proper analogue of the Gaussian for the
Askey-Wilson theory, and
we compute the image under the polynomial Askey-Wilson transform
of an Askey-Wilson polynomial multiplied by the inverse of the Gaussian.
In the special case of continuous
$q$-ultraspherical polynomials, these formulas were
derived by Cherednik in \cite{C2}.
Furthermore, we compute the image under the Askey-Wilson
function transform of an Askey-Wilson polynomial multiplied by the Gaussian.
A special case leads 
to the evaluation of a $q$-analogue (in one variable) of the
Macdonald-Mehta integral
(cf. Macdonald \cite{Macdonald}).

The techniques employed in this paper are entirely
based on basic hypergeometric series manipulations
in the spirit of Gasper's and Rahman's \cite{GR} book.
The two main ingredients are the orthogonality relations
for the Askey-Wilson polynomials (see \cite{AW}), and the inversion
formula for the Askey-Wilson function transform (see \cite{KS2}).

A generalization of Cherednik's \cite{C}
affine Hecke algebra approach to the Askey-Wilson level 
leads to independent proofs of the
Plancherel and inversion formula for the Askey-Wilson function
transform, and to independent proofs of
most of the formulas presented in this paper. In fact,
the affine Hecke algebra approach reduces the problem to
the explicit evaluation of the 
$q$-analogue of the (one variable) Macdonald-Mehta 
integral. I therefore have added 
two alternative proofs of the evaluation of the (one variable) 
$q$-Macdonald-Mehta integral in this paper, which do not make use of
the properties of the Askey-Wilson function transform. 
I will discuss the affine Hecke algebra approach in a future paper.

The plan of the paper is as follows.
In Section 2 we recall the basic properties
of the Askey-Wilson polynomials.
In Section 3 we give the definition of the Askey-Wilson
function. The expansion formula for the Askey-Wilson function
in terms of Askey-Wilson polynomials is formulated in Section 4.
In Section 4 we also introduce the analogue of the Gaussian, and
we explicitly compute the image under the polynomial Askey-Wilson
transform of an Askey-Wilson polynomial multiplied by the inverse of the
Gaussian. In Section 5 the Askey-Wilson function transform and its
basic properties are recalled,
and the image under the Askey-Wilson function transform
of an Askey-Wilson polynomial multiplied by the Gaussian is computed
explicitly. We also show how this leads to the evaluation of a
$q$-analogue (in one variable) of the Macdonald-Mehta integral.
In Section 6 some density results are discussed, 
which are relevant for the $L^2$-theory of
the Askey-Wilson function transform. This leads to
explicit parameter restraints for which
the formulas derived in Section 5 
completely determine the Askey-Wilson function transform.
Appendix A contains a proof of (a reformulation
of) the expansion formula for the Askey-Wilson function. 
Appendix B contains two direct proofs for the evaluation of
the $q$-analogue of the (one variable) Macdonald-Mehta integral.

%%%%%%%%%%%%%%%%%%%%%%%%%%%%%%%%%%%%%%%%%%%%%%%%%%%%%%%%%%%%%%%%
{\it Notations and conventions:}
Throughout the paper we fix $0<q<1$.
The notation $\mathbb{C}^\times$ and $\mathbb{R}^\times$
is used for $\mathbb{C}\setminus \{0\}$ and $\mathbb{R}\setminus \{0\}$,
respectively.
The non-negative integers $\{0,1,2,\ldots\}$ are
denoted by $\mathbb{Z}_+$. The book \cite{GR} of Gasper and Rahman
is used as main reference for notations and results concerning basic
hypergeometric series. 
For $k\in\mathbb{Z}_+\cup \{\infty\}$ we write
$\bigl(x_1,\ldots,x_r;q\bigr)_{k}=\prod_{i=1}^r\bigl(x_i;q\bigr)_k$
with $\bigl(x;q\bigr)_k=\prod_{i=0}^{k-1}(1-xq^i)$ the $q$-shifted
factorial. Similarly, we write
$\theta(a_1,\ldots,a_r)=\prod_{i=1}^r\theta(a_i)$ with
$\theta(a)=\bigl(a,q/a;q\bigr)_{\infty}$ for (products) of
renormalized Jacobi theta functions. The series expansion
\begin{equation*}
{}_{r}\phi_s\left(\begin{matrix} a_1, a_2,\ldots, a_{r}\\
b_1,b_2,\ldots,b_s \end{matrix}\,; q,z\right)=
\sum_{k=0}^\infty\frac{\bigl(a_1,a_2,\ldots,a_{r};q\bigr)_k}{\bigl(q,b_1,
\ldots,b_s;q\bigr)_k}\lbrack (-1)^kq^{\frac{1}{2}k(k-1)}\rbrack^{1+s-r}z^k
\end{equation*}
defines the ${}_{r}\phi_s$ basic hypergeometric series.
The very-well-poised ${}_8\phi_7$ basic hypergeometric
series is defined by
\begin{equation*}
{}_8W_7(a;b,c,d,e,f;q,z)=
\sum_{k=0}^{\infty}\frac{1-aq^{2k}}{1-a}
\frac{\bigl(a,b,c,d,e,f;q\bigr)_kz^k}{\bigl(q,qa/b,qa/c,qa/d,qa/e,
qa/f;q\bigr)_k}.
\end{equation*}
The bilateral basic hypergeometric series ${}_r\psi_s$ is defined by
\begin{equation*}
{}_r\psi_s\left(\begin{matrix} a_1, a_2,\ldots, a_{r}\\
b_1,b_2,\ldots,b_s \end{matrix}\,;q,z\right)=
\sum_{n\in \mathbb{Z}}\frac{\bigl(a_1,a_2,\ldots,a_r;q\bigr)_n}
{\bigl(b_1,b_2,\ldots,b_s;q\bigr)_n}\lbrack (-1)^nq^{\frac{1}{2}n(n-1)}
\rbrack^{s-r}\,z^n.
\end{equation*}
We use the branch of the square root $\sqrt{\,\cdot\,}$\,
which is positive on ${\mathbb{R}}_{>0}$, with branch cut along
the half-line $(-\infty,0)$ of the complex plane.

%%%%%%%%%%%%%%%%%%%%%%%%%%%%%%%%%%%%%%%%%%%%%%%%%%%%%%%%

{\it Acknowledgements:} The author is
supported by a fellowship from the Royal
Netherlands Academy of Arts and Sciences (KNAW).
I thank Ivan Cherednik, Marcel de Jeu, Erik Koelink and Eric Opdam for
stimulating discussions. 
%%%%%%%%%%%%%%%%%%%%%%%%%%%%%%%%%%%%%%%%%%%%%%%%%%%%%%%%%%%%%%%%%%%%%%

%%%%%%%%%%%%%%%%%%%%%%%%%%%%%%%%%%%%%%%%%%%
%%                                       %%
%%    The Askey-Wilson polynomials       %%
%%                                       %%
%%%%%%%%%%%%%%%%%%%%%%%%%%%%%%%%%%%%%%%%%%%

\section{The Askey-Wilson polynomials}

In order to fix notations, we recall
the basic properties of the Askey-Wilson
polynomials in this section.

The Askey-Wilson polynomials depend,
besides $q$, on four parameters $a,b,c,d$. To simplify notations
it is convenient to use the short-hand notation
\[ \alpha=(a,b,c,d)
\]
for the four-tuple of parameters $a,b,c,d$, which we assume throughout
this section to be generically complex and subject to the condition
$\hbox{Re}(a)>0$.
We define dual parameters
\begin{equation}\label{dualpar}
\alpha_\sigma=(a_\sigma,b_\sigma,c_\sigma,d_\sigma)
\end{equation}
by
\[ a_\sigma=\sqrt{q^{-1}abcd},\qquad
b_\sigma=ab/a_\sigma,\qquad c_\sigma=ac/a_\sigma,
\qquad d_\sigma=ad/a_\sigma.
\]
This notation turns out to be quite useful
later on when we have to compose involutions on parameter
sets. Since dual parameters play an important role throughout the paper, it
is convenient to have a second, less cumbersome notation
at hand. This second notation is
\[ \alpha_\sigma=\widetilde{\alpha},\qquad
(a_\sigma,b_\sigma,c_\sigma,d_\sigma)=
(\tilde{a},\tilde{b},\tilde{c},\tilde{d}),
\]
in accordance with \cite{KS2}. The map
$\alpha\mapsto \widetilde{\alpha}$
defines an involution on the four tuple of parameters $\alpha$.
Here the condition $\hbox{Re}(a)>0$ is needed
in view of the chosen branch for $\sqrt{\,\cdot\,}$,
see the conventions at the end of the introduction. Observe
in particular that
$\hbox{Re}(\tilde{a})=\hbox{Re}(\sqrt{q^{-1}abcd})>0$ for generic parameters
$\alpha$ in view of
the chosen branch for $\sqrt{\,\cdot\,}$.
\begin{rem}
Throughout the paper we
formulate the results under the assumption $\hbox{Re}(a)>0$ in order to be
able to use the duality involution $\sigma$ without worrying about the
chosen branch of the square-root. In most formulas the condition
$\hbox{Re}(a)>0$ can be easily removed by analytic continuation.
\end{rem}

We define a discrete subset
$\mathcal{S}=\mathcal{S}(\alpha;q)\subset \mathbb{C}^\times$ by
\[\mathcal{S}=\{ s_m \, | \, m\in \mathbb{Z}_+\},\qquad
s_m=s_m(\alpha;q)=\tilde{a}q^m.
\]

The Askey-Wilson polynomials
$E_s(x)=E_s(x;\alpha;q)$\, ($s\in\mathcal{S}$)
are defined by the series expansion
\begin{equation}\label{pn}
\begin{split}
E_{s_n}(x)&={}_4\phi_3\left(\begin{matrix} \tilde{a}\,s_n,
\tilde{a}/s_n, ax, a/x \\
ab, ac, ad \end{matrix}\,; q,q\right)\\
&={}_4\phi_3\left(\begin{matrix} q^{-n},  q^{n-1}abcd, ax, a/x \\
 ab, ac, ad \end{matrix}\,; q,q\right)
\end{split}
\end{equation}
for $n\in\mathbb{Z}_+$, see \cite{AW}.
For fixed $s\in\mathcal{S}$, the Askey-Wilson polynomial $E_s(x)$
is an eigenfunction of the Askey-Wilson second-order $q$-difference
operator $L=L(\alpha;q)$,
\begin{equation}\label{LAW}
\begin{split}
(Lp)(x)&=C(x)(p(qx)-p(x))+C(x^{-1})(p(q^{-1}x)-p(x)),\\
C(x)&=\frac{(1-ax)(1-bx)(1-cx)(1-dx)}{(1-x^2)(1-qx^2)},
\end{split}
\end{equation}
with eigenvalue $\mu(s)$, where
\[\mu(\gamma)=-1-\tilde{a}^2+\tilde{a}(\gamma+\gamma^{-1}).
\]
Furthermore, the Askey-Wilson polynomial $E_s(x)$
has the duality property
\begin{equation}\label{dualitypol}
E_s(v)=\widetilde{E}_v(s),\qquad s\in \mathcal{S},\,\,\, v\in
\widetilde{\mathcal{S}},
\end{equation}
where $\widetilde{\mathcal{S}}=\mathcal{S}(\widetilde{\alpha};q)$
and
$\widetilde{E}_v(\cdot)=E_v(\cdot;\widetilde{\alpha};q)$
for $v\in \widetilde{\mathcal{S}}$.

Let
$\mathcal{T}=\mathcal{T}_{\alpha,q}$ be a closed, counterclockwise oriented
contour in the complex plane, for which the sequences $eq^{\mathbb{Z}_+}$
(respectively $e^{-1}q^{-\mathbb{Z}_+}$) are in the interior
(respectively exterior) of $\mathcal{T}$ for all $e=a,b,c,d$.
In case $|a|,|b|,|c|,|d|<1$, one can for instance take for $\mathcal{T}$
the unit circle $\mathbb{T}$ in the complex plane.

We call a function $f:\mathbb{C}^\times\rightarrow \mathbb{C}$
{\it inversion-invariant} if $f(x^{-1})=f(x)$ for all
$x\in\mathbb{C}^\times$.
For ``sufficiently nice'' inversion-invariant functions $f$
we define the {\it polynomial Askey-Wilson transform}
$(\mathbb{F}f)(s)=(\mathbb{F}(\alpha,q)f)(s)$
of $f$ at $s\in \mathcal{S}$ by
\begin{equation*}
(\mathbb{F}f)(s)=
\frac{1}{4\pi i N}\int_\mathcal{T}f(x)E_s(x)\Delta(x)\frac{dx}{x},
\end{equation*}
where $\Delta(x)=\Delta(x;\alpha;q)$ is the
weight function
\[ \Delta(x)=\frac{\bigl(x^2,1/x^{2};q\bigr)_{\infty}}
{\bigl(ax,a/x,bx,b/x,cx,c/x,dx,d/x;q\bigr)_{\infty}},
\]
and the constant $N=N(\alpha;q)$ is the Askey-Wilson integral
\begin{equation}\label{N}
\begin{split}
N&=\frac{1}{4\pi i}\int_{\mathcal{T}}\Delta(x)\frac{dx}{x}\\
&=\frac{\bigl(abcd;q\bigr)_{\infty}}
{\bigl(q,ab,ac,ad,bc,bd,cd;q\bigr)_{\infty}},
\end{split}
\end{equation}
see \cite{AW}.
Conversely, for ``sufficiently nice'' functions
$g:\mathcal{S}\rightarrow \mathbb{C}$ we define the transform
$(\mathbb{I}g)(x)=(\mathbb{I}(\alpha,q)g)(x)$\, ($x\in\mathbb{C}^\times$) by
\[(\mathbb{I}g)(x)=\sum_{s\in \mathcal{S}}g(s)E_s(x)\widetilde{h}(s),
\]
with the weight $\widetilde{h}(s)=h(s;\widetilde{\alpha};q)$
for $s\in\mathcal{S}$ given by
\begin{equation}\label{hnorm}
\begin{split}
\widetilde{h}(s_m)&=\frac{\underset{x=s_m}{\textup{\hbox{Res}}}
\left(\frac{\widetilde{\Delta}(x)}{x}\right)}
{\underset{x=s_0}{\textup{\hbox{Res}}}
\left(\frac{\widetilde{\Delta}(x)}{x}\right)}\\
&=\frac{\bigl(1-q^{2m-1}abcd\bigr)\bigl(q^{-1}abcd,ab,
ac,ad;q\bigr)_{m}}{\bigl(1-q^{-1}abcd\bigr)\bigl(q,bc,bd,cd;q\bigr)_m}a^{-2m}
\end{split}
\end{equation}
for $m\in \mathbb{Z}_+$, cf. \cite{NS}. Here
$\widetilde{\Delta}(x)=\Delta(x;\widetilde{\alpha};q)$
is the weight function $\Delta(x)$ with respect to dual parameters.
In this paper, we consider the transform $\mathbb{F}$ respectively
$\mathbb{I}$ with respect to two classes of functions $f$ respectively $g$.
We first consider the
function space $\mathcal{A}=\mathbb{C}[x+x^{-1}]$
consisting of inversion-invariant Laurent polynomials in the
variable $x$ for $\mathbb{F}$. Observe that the Askey-Wilson polynomials
$\{E_s \, | \, s\in\mathcal{S} \}$ form a linear basis of $\mathcal{A}$.
The corresponding function space
$\mathcal{F}_0(\mathcal{S})$ for $\mathbb{I}$ consists of
functions $g: \mathcal{S}\rightarrow \mathbb{C}$
with finite support. The set of delta-functions $\{ \delta_s \, |
\, s\in\mathcal{S} \}$, with $\delta_s(v)=\delta_{s,v}$
for $s,v\in \mathcal{S}$, forms a linear basis of
$\mathcal{F}_0(\mathcal{S})$.
The orthogonality relations
\begin{equation}\label{orthogonality}
\frac{1}{4\pi i
N}\int_{\mathcal{T}}E_s(x)E_v(x)\Delta(x)\frac{dx}{x}=
\delta_{s,v}\,\frac{1}{\widetilde{h}(s)},\qquad s,v\in
\mathcal{S},
\end{equation}
for the Askey-Wilson polynomials (see \cite[Theorem 2.3]{AW}),
imply that $\mathbb{F}(E_s)=\widetilde{h}(s)^{-1}\delta_s$ for
$s\in\mathcal{S}$. On the other hand, $\mathbb{I}(\delta_s)=
\widetilde{h}(s)E_s$ by the definition of $\mathbb{I}$.
This immediately leads to the following theorem.
\begin{thm}\label{AWresult}
$\mathbb{F}$ defines a linear bijection
$\mathbb{F}:\mathcal{A}\rightarrow \mathcal{F}_0(\mathcal{S})$. Its
inverse is given by
$\mathbb{I}: \mathcal{F}_0(\mathcal{S})\rightarrow \mathcal{A}$.
\end{thm}

%%%%%%%%%%%%%%%%%%%%%%%%%%%%%%%%%%%%%%%%%%%%%%%%
%%                                            %%
%%   The Askey-Wilson function                %%
%%                                            %%
%%%%%%%%%%%%%%%%%%%%%%%%%%%%%%%%%%%%%%%%%%%%%%%%

\section{The Askey-Wilson function}
In this section we recall the definition of the
Askey-Wilson function (see e.g. \cite{S89}, \cite{S},\cite{S2}, \cite{IR}
and \cite{KS2}), and give some of its basic properties. 
The Askey-Wilson function is a non-polynomial eigenfunction of
the Askey-Wilson second order $q$-difference operator, given
explicitly as a very-well-poised ${}_8\phi_7$ series.  
Recall that an explicit basis of eigenfunctions for the
Askey-Wilson second-order $q$-difference operator
in terms of very-well-poised  ${}_8\phi_7$ series is known,
see Ismail and Rahman \cite{IR} (compare also with Suslov \cite{S89}).
See Ruijsenaars \cite{R} for the case $|q|=1$, which requires a completely
different approach.

We assume
that the parameters $\alpha=(a,b,c,d)$ are generically complex, and
subject to the condition $\hbox{Re}(a)>0$.
The Askey-Wilson function
$\phi_\gamma(x)=\phi_\gamma(x;\alpha;q)$ is defined by
\begin{equation}\label{phi}
\begin{split}
\phi_{\gamma}(x)=&\frac{\bigl(qax\gamma/\tilde{d},
qa\gamma/\tilde{d}x,qabc/d;q\bigr)_{\infty}}
{\bigl(\tilde{a}\tilde{b}\tilde{c}\gamma, q\gamma/\tilde{d},
qx/d, q/dx,bc,qb/d,qc/d;q\bigr)_{\infty}}\\
&\qquad\times{}_8W_7\bigl(\tilde{a}\tilde{b}\tilde{c}\gamma/q;
ax, a/x, \tilde{a}\gamma,
\tilde{b}\gamma, \tilde{c}\gamma\,; q,q/\tilde{d}\gamma\bigr),\qquad
|q/\tilde{d}\gamma|<1.
\end{split}
\end{equation}
Note that $\phi_\gamma(x)$ is normalized differently
compared with \cite{KS2}. It is known that
\[(L\phi_\gamma)(x)=\mu(\gamma)\phi_\gamma(x),
\]
where $L$ is the Askey-Wilson second-order
$q$-difference operator \eqref{LAW}, see e.g. \cite{IR}, \cite{S}, \cite{KS2}.
In view of Bailey's formula \cite[(2.10.10)]{GR} we can write
\begin{equation}\label{43presentation}
\begin{split}
\phi_{\gamma}(x)=&\frac{\bigl(qabc/d;q\bigr)_{\infty}}
{\bigl(bc,qa/d,qb/d,qc/d,q/ad;q\bigr)_{\infty}}
{}_4\phi_3\left(\begin{matrix} ax, a/x, \tilde{a}\gamma, \tilde{a}/\gamma\\
 ab, ac, ad \end{matrix}\,; q,q\right)\\
+&\frac{\bigl(ax, a/x, \tilde{a}\gamma, \tilde{a}/\gamma, qabc/d
;q\bigr)_{\infty}}
{\bigl(qx/d, q/dx, q\gamma/\tilde{d}, q/\tilde{d}\gamma,
ab,ac,bc,qa/d,ad/q;q\bigr)_{\infty}}\\
&\qquad\qquad\qquad\quad\times
{}_4\phi_3\left(\begin{matrix} qx/d, q/dx, q\gamma/\tilde{d},
q/\tilde{d}\gamma\\
 qb/d, qc/d, q^2/ad \end{matrix}\,; q,q\right).
\end{split}
\end{equation}
In particular,
$\phi_{\gamma}(x)$ extends to a meromorphic function in
$(\gamma,x)\in\mathbb{C}^\times\times\mathbb{C}^\times$
and is inversion-invariant
in both $x$ and $\gamma$.
The possible poles of $\phi_\gamma(x)$ are simple and are located at
$x^{\pm 1}=q^{1+k}/d$\, ($k\in {\mathbb{Z}}_+$) and
$\gamma^{\pm 1}=q^{1+k}/\tilde{d}$\, ($k\in {\mathbb{Z}}_+$).
It follows from \eqref{43presentation} that
\begin{equation}\label{redpol}
 \phi_{s}(x)=\frac{\bigl(qabc/d;q\bigr)_{\infty}}
{\bigl(bc,q/ad,qa/d,qb/d,qc/d;q\bigr)_{\infty}}E_s(x),
\qquad s\in \mathcal{S},
\end{equation}
and that
\begin{equation}\label{dualityfun}
\phi_\gamma(x)=\widetilde{\phi}_x(\gamma)
\end{equation}
where $\widetilde{\phi}_x(\gamma)=\phi_x(\gamma;\widetilde{\alpha};q)$
is the Askey-Wilson function with respect to dual parameters.
Formula \eqref{redpol} implies that the Askey-Wilson function is a
meromorphic continuation of the Askey-Wilson polynomial in its
degree. We will refer to formula
\eqref{redpol} as the {\it polynomial reduction}
of the Askey-Wilson function. Formula \eqref{dualityfun} implies
that the geometric parameter $x$ and the spectral parameter $\gamma$ of the
Askey-Wilson function are interchangeable in a suitable sense. We will refer 
to formula \eqref{dualityfun} as the
{\it duality} of the Askey-Wilson function.
It extends the duality \eqref{dualitypol} of the Askey-Wilson polynomial.

%%%%%%%%%%%%%%%%%%%%%%%%%%%%%%%%%%%%%%%%%%%%%%%%%%%%%%%%%%%%%%%%%
%%                                                             %%
%%       The expansion formula                                 %%
%%                                                             %%
%%%%%%%%%%%%%%%%%%%%%%%%%%%%%%%%%%%%%%%%%%%%%%%%%%%%%%%%%%%%%%%%%

\section{The expansion formula.}

We assume in this section that the parameters
$\alpha=(a,b,c,d)$ are generically complex and subject to the
condition $\hbox{Re}(a)>0$.
In order to formulate the expansion formula for the
Askey-Wilson function
in terms of Askey-Wilson polynomials, it is important to
keep track of two involutions on the four tuples
$\alpha=(a,b,c,d)$. Firstly we have the concept of dual parameters,
which we have
already used extensively in the previous sections. It is now more
convenient to write the dual parameter with sub-index $\sigma$, so
\[ \alpha_\sigma=(a_\sigma,b_\sigma,c_\sigma,d_\sigma)=
(\tilde{a},\tilde{b},\tilde{c},\tilde{d})=\widetilde{\alpha},
\]
with $\alpha_\sigma$ defined by \eqref{dualpar}.
Secondly, we define $\tau$ by
\begin{equation}
\alpha_\tau=(a_\tau,b_\tau,c_\tau,d_\tau)=(a,b,c,q/d).
\end{equation}
We admit compositions of $\sigma$ and $\tau$, for instance
we write
\[
\alpha_{\sigma\tau}=
(a_{\sigma\tau},b_{\sigma\tau},c_{\sigma\tau},d_{\sigma\tau})
\]
for first applying $\sigma$ to $\alpha$,
and then applying $\tau$ to $\alpha_\sigma$, i.e.
\[\alpha_{\sigma\tau}=(\tilde{a},\tilde{b},\tilde{c},q/\tilde{d}).
\]
Observe that
\begin{equation}\label{parameterinvariant}
\alpha_{\sigma\tau\sigma}=\alpha_{\tau\sigma\tau},
\qquad \alpha_{\tau\tau}=\alpha.
\end{equation}
Furthermore, we have seen in Section 2 that
$\alpha_{\sigma\sigma}=\alpha$ since $\hbox{Re}(a)>0$.

Finally we use the convention
that if $H=H(\alpha)$ is an object depending on $\alpha$,
then e.g. $H^{\sigma\tau}$, or $H_{\sigma\tau}$, denotes the same object
in which the four tuple $\alpha$ is replaced by $\alpha_{\sigma\tau}$.
We sometimes write $\widetilde{H}=H^\sigma$ to simplify notations.

We define the {\it Gaussian} $G(x)=G(x;\alpha;q)$ by
\begin{equation}
G(x)=\bigl(dx,d/x;q\bigr)_{\infty}^{-1}.
\end{equation}
The terminology stems from Cherednik's \cite{C}, \cite{C2}
work on Gaussians associated with Macdonald polynomials. 

The {\it analytic part}
$\phi_\gamma^{an}(x)=\phi_\gamma^{an}(x;\alpha;q)$
of the Askey-Wilson function
$\phi_\gamma(x)=\phi_\gamma(x;\alpha;q)$ is defined by
\begin{equation}\label{holomorphicpart}
\begin{split}
\phi_\gamma^{an}(x)&=
G^\tau(x)^{-1}G^{\sigma\tau}(\gamma)^{-1}\phi_{\gamma}(x)\\
&=\bigl(qx/d,q/dx,q\gamma/\tilde{d},q/\tilde{d}\gamma;q\bigr)_{\infty}\,
\phi_\gamma(x).
\end{split}
\end{equation}
The properties of $\phi_\gamma(x)$ as described in Section 3 imply that
$\phi_\gamma^{an}(x)$ is
analytic in $(\gamma,x)\in \mathbb{C}^\times\times
\mathbb{C}^\times$.

Finally we observe that $\mathcal{S}^\tau=\mathcal{S}^{\sigma\tau}$, since
\begin{equation}\label{sst}
s_m^\tau=s_m^{\sigma\tau}=q^m\sqrt{abc/d},\qquad m\in\mathbb{Z}_+.
\end{equation}
We can now formulate the following key proposition.
\begin{prop}\label{expansion2}
For $s\in\mathcal{S}^\tau=\mathcal{S}^{\sigma\tau}$ and
$\gamma\in\mathbb{C}^\times$,
\begin{equation}\label{secondform}
\bigl(\mathbb{F}^\tau(\phi_\gamma^{an})\bigr)(s)=
\frac{G^{\tau\sigma\tau}(s)}{G^{\tau\sigma\tau}(s_0^\tau)}
E_s^{\sigma\tau}(\gamma).
\end{equation}
\end{prop}
The proof of the proposition, which is based on direct calculations
using the theory of basic hypergeometric series, is given in Appendix A.
Proposition \ref{expansion2} leads to the following
expansion formula for the analytic part $\phi_\gamma^{an}(x)$
of the Askey-Wilson function.
\begin{thm}[The expansion formula]\label{expansion}
\begin{equation}\label{equivform}
\begin{split}
\phi_\gamma^{an}(x)&=
G^{\tau\sigma\tau}(s_0^\tau)^{-1}\sum_{s\in \mathcal{S}^\tau}
h^{\tau\sigma}(s)G^{\tau\sigma\tau}(s)
E^{\sigma\tau}_s(\gamma)E^\tau_s(x)\\
&=\sum_{m=0}^\infty
{}_{4}\phi_3\left(\begin{matrix} q^{-m}, q^mabc/d, ax, a/x\\
ab,ac,qa/d\end{matrix}\,; q,q\right)
{}_{4}\phi_3\left(\begin{matrix} q^{-m},
q^mabc/d, \tilde{a}\gamma,
\tilde{a}/\gamma\\
ab,ac,bc \end{matrix}\,; q,q\right)\\
&\qquad\qquad\qquad\times \frac{(1-abcq^{2m}/d)\bigl(abc/d,ab,ac;q\bigr)_m}
{(1-abc/d)\bigl(q,qb/d,qc/d;q\bigr)_m}\left(\frac{-1}{ad}\right)^mq^{m(m+1)/2}
\end{split}
\end{equation}
for all $(\gamma,x)\in\mathbb{C}^\times\times\mathbb{C}^\times$.
\end{thm}
\begin{proof}
First observe that the terms 
$h^{\tau\sigma}(s)E^{\sigma\tau}_s(\gamma)E^\tau_s(x)$
for $s\in\mathcal{S}^\tau$ in the expansion sum
\eqref{equivform} are well defined in view of \eqref{sst}.
Observe furthermore that the second
expansion sum in \eqref{equivform} converges absolutely
and uniformly for $(\gamma,x)$ in compacta of
$\mathbb{C}^\times\times\mathbb{C}^\times$ due to the
Gaussian type factor $q^{m(m+1)/2}$
(use e.g. \cite[(7.5.13)]{GR} to control the convergence of the
${}_4\phi_3$'s in the expansion sum). In particular, the second
expansion sum in \eqref{equivform} is analytic in
$(\gamma,x)\in\mathbb{C}^\times\times\mathbb{C}^\times$.
We can now verify the second identity in \eqref{equivform} term-wise using
the explicit expression
for the Askey-Wilson polynomial as a balanced ${}_4\phi_3$ series
(see \eqref{pn}), and using the identities
\[\frac{G^{\tau\sigma\tau}(s_m^\tau)}
{G^{\tau\sigma\tau}(s_0^\tau)}=
\frac{\bigl(bc;q\bigr)_m}
{\bigl(qa/d;q\bigr)_m}\left(\frac{-a}{d}\right)^mq^{m(m+1)/2}
\]
and
\[ h^{\tau\sigma}(s_m^\tau)=\frac{\bigl(1-abcq^{2m}/d\bigr)
\bigl(abc/d,ab,ac,qa/d;q\bigr)_m}
{\bigl(1-abc/d\bigr)\bigl(q,bc,qb/d,qc/d;q\bigr)_m}a^{-2m}
\]
for $m\in\mathbb{Z}_+$.
So it remains to prove the first identity in \eqref{equivform}.
Denote $\psi_\gamma(x)$ for the right hand side of
\eqref{equivform}, which we consider for arbitrary, fixed
$\gamma\in\mathbb{C}^\times$ as an analytic, inversion-invariant
function in $x\in\mathbb{C}^\times$. Recall that the defining expansion sum
for $\psi_\gamma(x)$ converges absolutely and uniformly for $x$ in
compacta of $\mathbb{C}^\times$. In particular, when applying the
polynomial Askey-Wilson transform $\mathbb{F}^\tau$ to $\psi_\gamma$,
it is allowed to interchange summation and integration. Combined
with the orthogonality relations \eqref{orthogonality} for the Askey-Wilson
polynomials, we obtain for $s\in \mathcal{S}^\tau$,
\begin{equation*}
\begin{split}
\bigl(\mathbb{F}^\tau(\psi_\gamma)\bigr)(s)&=
G^{\tau\sigma\tau}(s_0^\tau)^{-1}\sum_{v\in \mathcal{S}^\tau}
h^{\tau\sigma}(v)G^{\tau\sigma\tau}(v)E_v^{\sigma\tau}(\gamma)
\left(\mathbb{F}^\tau(E_v^\tau)\right)(s)\\
&=\frac{G^{\tau\sigma\tau}(s)}{G^{\tau\sigma\tau}(s_0^\tau)}\,
E_s^{\sigma\tau}(\gamma)\\
&=\bigl(\mathbb{F}^\tau(\phi_\gamma^{an})\bigr)(s)
\end{split}
\end{equation*}
where the last equality follows from
Proposition \ref{expansion2}. Since any analytic, inversion-invariant
function $f:\mathbb{C}^\times\rightarrow\mathbb{C}$ is
uniquely determined by its image under the polynomial Askey-Wilson
transform $\mathbb{F}^\tau$, we conclude that
$\psi_\gamma(x)=\phi_\gamma^{an}(x)$ for all $x\in
\mathbb{C}^\times$,
as desired.
\end{proof}

\begin{rem}
The right hand side of the expansion formula \eqref{equivform}
resembles the (non-symmetric) Poisson-kernel
for Askey-Wilson polynomials, see \cite{ARS}.
The essential difference is the occurrence
of the Gaussian $G^{\tau\sigma\tau}(s)$ in the
expansion sum \eqref{equivform}, which is not present in the Poisson-kernels.
It is also due to this extra factor that the expansion sum \eqref{equivform}
has better convergence properties.
\end{rem}

We explore now the implications of the formulas \eqref{secondform} for the
polynomial Askey-Wilson transform $\mathbb{F}$ and its inverse
$\mathbb{I}$. Let $\mathcal{F}(\mathcal{S})$ be the space of
functions $g: \mathcal{S}\rightarrow \mathbb{C}$.
\begin{prop}\label{compactcontmaster}
For $s\in \mathcal{S}^\tau=\mathcal{S}^{\sigma\tau}$, we have
\begin{equation}\label{compactcontmasterform}
\mathbb{F}\bigl(E_s^{\tau}G^{-1}\bigr)=
\frac{\bigl(bc,bc,d/a,ad,q/ad,bd,cd;q\bigr)_{\infty}}
{\bigl(abcd;q\bigr)_{\infty}}\,G^{\tau\sigma\tau}(s)
E_s^{\sigma\tau}G^{\sigma\tau}\in \mathcal{F}(\mathcal{S}).
\end{equation}
\end{prop}
\begin{proof}
Let $s\in \mathcal{S}^\tau$ and $v\in \mathcal{S}$.
It follows from the polynomial reduction \eqref{redpol} of the
Askey-Wilson function, \eqref{holomorphicpart} and \eqref{secondform},
that
\begin{equation*}
\begin{split}
\bigl(\mathbb{F}\bigl(E_s^{\tau}G^{-1}\bigr)\bigr)(v)&=
\frac{\bigl(bc,d/a,ad,bd,cd;q\bigr)_{\infty}}
{\bigl(abcd;q\bigr)_{\infty}}\,
\bigl(\mathbb{F}\bigl(\phi_s^\tau G^{-1}\bigr)\bigr)(v)\\
&=\frac{\bigl(bc,d/a,ad,bd,cd;q\bigr)_{\infty}}
{\bigl(abcd;q\bigr)_{\infty}}\,
G^{\tau\sigma\tau}(s)\bigl(\mathbb{F}\bigl(\phi_s^{an,\tau}\bigr)\bigr)(v)\\
&=\frac{\bigl(bc,d/a,ad,bd,cd;q\bigr)_{\infty}}
{\bigl(abcd;q\bigr)_{\infty}}\,\frac{G^{\tau\sigma\tau}(s)G^{\sigma\tau}(v)}
{G^{\sigma\tau}(s_0)}E_v^{\tau\sigma\tau}(s).
\end{split}
\end{equation*}
Now $\mathcal{S}=\mathcal{S}^{\sigma\tau\sigma}$
and $\mathcal{S}^\tau=\mathcal{S}^{\sigma\tau}=
\mathcal{S}^{\sigma\tau\sigma\sigma}$, so we
conclude from \eqref{parameterinvariant} and from
the duality \eqref{dualitypol} of the Askey-Wilson
polynomials that
\begin{equation}\label{sideremark}
E_v^{\tau\sigma\tau}(s)=E_v^{\sigma\tau\sigma}(s)=E_s^{\sigma\tau}(v).
\end{equation}
Furthermore, $G^{\sigma\tau}(s_0)=\bigl(bc,q/ad;q\bigr)_{\infty}^{-1}$,
hence
\[\bigl(\mathbb{F}(E_s^{\tau}G^{-1})\bigr)(v)
=\frac{\bigl(bc,bc,d/a,ad,q/ad,bd,cd;q\bigr)_{\infty}}
{\bigl(abcd;q\bigr)_{\infty}}\,G^{\tau\sigma\tau}(s)
E_s^{\sigma\tau}(v)G^{\sigma\tau}(v),
\]
as desired.
\end{proof}
\begin{rem}\label{Hahn}
Observe that
\[ G(x)^{-1}\Delta(x)=\frac{\bigl(x^2,1/x^{2};q\bigr)_{\infty}}
{\bigl(ax,a/x,bx,b/x,cx,c/x;q\bigr)_{\infty}}
\]
is the orthogonality 
weight function for the continuous dual $q$-Hahn polynomials.
Formula \eqref{compactcontmasterform}
thus computes the constant term of the product of two Askey-Wilson
polynomials (with parameters $(a,b,c,q/d)$ and $(a,b,c,d)$ respectively)
with respect to the continuous dual $q$-Hahn
orthogonality measure. This observation leads to the following alternative
way to prove Proposition \ref{compactcontmaster}.
First use the explicit expansion of the Askey-Wilson polynomial
as linear combination of continuous dual $q$-Hahn polynomials
(see \cite[(7.6.8) \& (7.6.9)]{GR}, and be aware of the fact
that a factor $\bigl(q;q\bigr)_n$ is missing
in the numerator of \cite[(7.6.9)]{GR}),
and substitute these for the two Askey-Wilson polynomials
$E_s^\tau$\, ($s\in \mathcal{S}^\tau$) and $E_v$\, ($v\in \mathcal{S}$)
in $\bigl(\mathbb{F}(E_s^{\tau}G^{-1})\bigr)(v)$.
Using the orthogonality relations
for the continuous dual $q$-Hahn polynomials we arrive at a single
sum, which is easily seen to give the same
result as \eqref{compactcontmasterform}.
\end{rem}
Cherednik's formulas \cite[(1.15)]{C2} involving continuous
$q$-ultraspherical polynomials can now be generalized to the
level of Askey-Wilson polynomials as follows.
\begin{thm}\label{polmaster}
The polynomial Askey-Wilson transform
$\mathbb{F}=\mathbb{F}(\alpha;q)$
defines a linear bijection
\[
\mathbb{F}: \mathcal{A}G^{-1}\rightarrow
\bigl(\mathcal{A}G^{\sigma\tau}\bigr)|_\mathcal{S}
\subset \mathcal{F}(\mathcal{S}),
\]
with inverse
\[
\mathbb{I}: \bigl(\mathcal{A}G^{\sigma\tau}\bigr)|_\mathcal{S}\rightarrow
\mathcal{A}G^{-1}.\]
Explicitly, we have
\begin{equation}\label{polmasterformulas}
\begin{split}
\mathbb{F}\bigl(E_s^{\tau}G^{-1}\bigr)&=
\frac{\bigl(bc,bc,d/a,ad,q/ad,bd,cd;q\bigr)_{\infty}}
{\bigl(abcd;q\bigr)_{\infty}}\,G^{\tau\sigma\tau}(s)
E_s^{\sigma\tau}G^{\sigma\tau},\\
\mathbb{I}\bigl(E_s^{\sigma\tau}G^{\sigma\tau}\bigr)&=
\frac{\bigl(abcd;q\bigr)_{\infty}}
{\bigl(bc,bc,d/a,ad,q/ad,bd,cd;q\bigr)_{\infty}}\,
G^{\tau\sigma\tau}(s)^{-1}E_s^{\tau}G^{-1},
\end{split}
\end{equation}
for $s\in \mathcal{S}^\tau=\mathcal{S}^{\sigma\tau}$.
\end{thm}
\begin{proof}
In view of the previous proposition,
it suffices to prove the explicit formula for
$\mathbb{I}\bigl(E_s^{\sigma\tau}G^{\sigma\tau}\bigr)$.

Let $s\in \mathcal{S}^\tau=\mathcal{S}^{\sigma\tau}$,
then $\mathbb{I}\bigl(E_s^{\sigma\tau}G^{\sigma\tau}\bigr)(x)$
is given by a series expansion in Askey-Wilson polynomials
$E_v(x)$\, ($v\in \mathcal{S}$) which converges absolutely and
uniformly on compacta of $x\in\mathbb{C}^\times$, compare with the proof
of Theorem \ref{expansion}. In particular,
$\mathbb{I}\bigl(E_s^{\sigma\tau}G^{\sigma\tau}\bigr)$ is an
inversion-invariant, analytic function. Furthermore, when applying
$\mathbb{F}$ to
$\mathbb{I}\bigl(E_s^{\sigma\tau}G^{\sigma\tau}\bigr)$ we may
interchange summation and integration. The orthogonality relations
\eqref{orthogonality} for the Askey-Wilson polynomials then show
that
\[\mathbb{F}
\bigl(\mathbb{I}\bigl(E_s^{\sigma\tau}G^{\sigma\tau}\bigr)\bigr)=
E_s^{\sigma\tau}G^{\sigma\tau}.
\]
On the other hand,
\eqref{compactcontmasterform} shows
that $E_s^{\sigma\tau}G^{\sigma\tau}\in \mathcal{F}(\mathcal{S})$
is the image under $\mathbb{F}$ of the analytic, inversion-invariant function
\[\frac{\bigl(abcd;q\bigr)_{\infty}}
{\bigl(bc,bc,d/a,ad,q/ad,bd,cd;q\bigr)_{\infty}}\,
G^{\tau\sigma\tau}(s)^{-1}E_s^{\tau}G^{-1}.
\]
Since any inversion-invariant, analytic function 
$f:\mathbb{C}^\times \rightarrow \mathbb{C}$ is uniquely determined by
its image under the polynomial Askey-Wilson transform $\mathbb{F}$,
we conclude that
\[\mathbb{I}\bigl(E_s^{\sigma\tau}G^{\sigma\tau}\bigr)=
\frac{\bigl(abcd;q\bigr)_{\infty}}
{\bigl(bc,bc,d/a,ad,q/ad,bd,cd;q\bigr)_{\infty}}\,G^{\tau\sigma\tau}(s)^{-1}
E_s^{\tau}G^{-1},
\]
as desired.
\end{proof}

\begin{rem}
The explicit formula \eqref{polmasterformulas} for
$\mathbb{I}\bigl(E_s^{\sigma\tau}G^{\sigma\tau}\bigr)$\,
($s\in\mathcal{S}^{\sigma\tau}$) can also
derived from the polynomial reduction \eqref{redpol} for the
Askey-Wilson function and from the
expansion formula \eqref{equivform}, since
\begin{equation*}
\begin{split}
E_s^{\tau}G^{-1}&=\frac{\bigl(bc,d/a,ad,bd,cd;q\bigr)_{\infty}}
{\bigl(abcd;q\bigr)_{\infty}}\,G^{-1} \phi_s^\tau \\
&=\frac{\bigl(bc,d/a,ad,bd,cd;q\bigr)_{\infty}}
{\bigl(abcd;q\bigr)_{\infty}}\,G^{\tau\sigma\tau}(s)\phi_s^{an,\tau}\\
&=\frac{\bigl(bc,d/a,ad,bd,cd;q\bigr)_{\infty}}
{\bigl(abcd;q\bigr)_{\infty}G^{\sigma\tau}(s_0)}\,G^{\tau\sigma\tau}(s)
\sum_{v\in\mathcal{S}}h^\sigma(v)G^{\sigma\tau}(v)E_v^{\tau\sigma\tau}(s)
E_v\\
&=\frac{\bigl(bc,bc,d/a,ad,q/ad,bd,cd;q\bigr)_{\infty}}
{\bigl(abcd;q\bigr)_{\infty}}\,G^{\tau\sigma\tau}(s)
\sum_{v\in\mathcal{S}}h^\sigma(v)G^{\sigma\tau}(v)E_s^{\sigma\tau}(v)
E_v\\
&=\frac{\bigl(bc,bc,d/a,ad,q/ad,bd,cd;q\bigr)_{\infty}}
{\bigl(abcd;q\bigr)_{\infty}}\,G^{\tau\sigma\tau}(s)\,
\mathbb{I}\bigl(E_s^{\sigma\tau}G^{\sigma\tau}\bigr),
\end{split}
\end{equation*}
where the fourth equality follows from \eqref{sideremark}.
\end{rem}

\begin{rem}\label{polMM}
The special case $s=s_0^\tau$ in the formula \eqref{polmasterformulas} for
$\mathbb{I}\bigl(E_s^{\sigma\tau}G^{\sigma\tau}\bigr)$ gives
\begin{equation*}
\begin{split}
&\bigl(dx,d/x;q\bigr)_{\infty}=
\frac{\bigl(bc,ad,q/ad,bd,cd;q\bigr)_{\infty}}
{\bigl(abcd;q\bigr)_{\infty}}\,
\bigl(\mathbb{I}(G^{\sigma\tau})\bigr)(x)\\
&=\frac{\bigl(ad,bd,cd;q\bigr)_{\infty}}
{\bigl(abcd;q\bigr)_{\infty}}\,
\sum_{m=0}^\infty
\frac{(1-q^{2m-1}abcd)\bigl(q^{-1}abcd,ab,ac;q\bigr)_m}
{(1-q^{-1}abcd)\bigl(q,bd,cd;q\bigr)_m}
\Bigl(\frac{-d}{q^{\frac{1}{2}}a}\Bigr)^mq^{\frac{m^2}{2}}E_{s_m}(x),
\end{split}
\end{equation*}
which may be viewed as a generalization of Jacobi's triple product identity
\cite[(1.6.1)]{GR} to the level of the Askey-Wilson
polynomials. Specializing $x=s_0^\sigma=a$ in this formula leads to
a limiting case of Rogers' \cite[(2.7.1)]{GR} summation formula of a
very-well-poised ${}_6\phi_5$ series,
\[{}_5\phi_5\left(\begin{matrix} \tilde{a}^2, q\tilde{a}, -q\tilde{a},
ab, ac\\
\tilde{a}, -\tilde{a}, bd, cd, 0\end{matrix}\,; q, d/a\right)
=\frac{\bigl(abcd,d/a;q\bigr)_{\infty}}{\bigl(bd,cd;q\bigr)_{\infty}}.
\]
The special case $\bigl(\mathbb{F}(G^{-1})\bigr)(\tilde{a})$ of
\eqref{polmasterformulas} is the evaluation of the
Askey-Wilson integral with one of the
four parameters equal to zero, cf. Remark \ref{Hahn}.
In Cherednik's \cite{C} terminology, both
$\bigl(\mathbb{F}(G^{-1})\bigr)(\tilde{a})$ and
$\bigl(\mathbb{I}(G^{\sigma\tau})\bigr)(a)$
are (polynomial) $q$-analogues of the (one variable)
Macdonald-Mehta integral.
\end{rem}

%%%%%%%%%%%%%%%%%%%%%%%%%%%%%%%%%%%%%%%%%%%%%%%%%%%%%%%%%%%%%%%%
%%                                                            %%
%%             The Askey-Wilson function transform            %%
%%                                                            %%
%%%%%%%%%%%%%%%%%%%%%%%%%%%%%%%%%%%%%%%%%%%%%%%%%%%%%%%%%%%%%%%%

\section{The Askey-Wilson function transform}

In \cite{KS2}, Koelink and the author
defined and studied a generalized Fourier transform
called the Askey-Wilson function transform, whose kernel is given
by the Askey-Wilson function.
In this section we show that Proposition \ref{expansion2},
together with the inversion formula \cite[Theorem 1]{KS2}
for the Askey-Wilson function transform, leads to an
explicit expression for the image under the Askey-Wilson function transform
of an Askey-Wilson polynomial multiplied by a Gaussian. These explicit
formulas lead
to a non-polynomial analogue of Theorem \ref{polmaster}.

We start by recalling the
definition and the main properties of the Askey-Wilson function
transform, see \cite{KS2} for more details. We use slightly different
conventions and normalizations compared with \cite{KS2}.

We fix a five tuple
\[\beta=(\alpha,t)=(a,b,c,d,t)\in {\mathbb{R}}^{5}\]
satisfying the conditions
\begin{equation}\label{conditions}
\begin{split}
&t<0,\qquad\qquad\qquad\qquad 0<b,c\leq a<d/q,\\
&bd,cd\geq q,\qquad\qquad\qquad ab,ac<1.
\end{split}
\end{equation}
The parameters $b,c,d$ then automatically satisfy $b,c<1$ and $d>q$.
The dual parameter $t_\sigma=\tilde{t}$ is defined by
\begin{equation}\label{dualt}
t_\sigma=\tilde{t}=\frac{1}{adt}.
\end{equation}
Note here the slightly different convention compared with
\cite[(4.4)]{KS2}.
The dual parameters $\beta_\sigma=\widetilde{\beta}=
(\widetilde{\alpha},\tilde{t})$
satisfies the same conditions \eqref{conditions} as
$\beta=(\alpha,t)$, see \cite[Lemma 1]{KS2}. In fact,
$\beta\mapsto \beta_\sigma$ defines an involution on the five
tuples $\beta=(\alpha,t)$ satisfying \eqref{conditions}.
For an object $H=H(\beta)$ depending on $\beta$, we write
$H^\sigma$, or $\widetilde{H}$, for $H(\beta_\sigma)$.

A new weight function
$W(x)=W(x;\beta;q)$ is defined by
\[ W(x)=\Delta(x)\Theta(x),
\]
where $\Delta(x)=\Delta(x;\alpha;q)$ is the weight function for
the Askey-Wilson polynomials and
$\Theta(x)=\Theta(x;\beta;q)$ is the quasi-constant
\begin{equation}\label{qcons}
\Theta(x)=\frac{\theta(dx,d/x)}{\theta(dtx,dt/x)}.
\end{equation}
For generic parameters $\beta$ such that the weight function $W$ 
has simple poles, we define a measure $m=m(\cdot\,;\beta;q)$ by
\begin{equation}\label{nu}
\int f(x)dm(x)=\frac{1}{4\pi i}\int_{x\in \mathbb{T}}f(x)
W(x)\frac{dx}{x}
+\frac{1}{2}\sum_{x\in D}\bigl(f(x)+f(x^{-1})\bigr)\,\underset{y=x}{\hbox{Res}}
\left(\frac{W(y)}{y}\right),
\end{equation}
where $D=D_+\cup D_-$ is the infinite, discrete set given by
\begin{equation}\label{D}
\begin{split}
D_+&=\{aq^k \, | \, k\in \mathbb{Z}_+: \, aq^k>1 \},\\
D_-&=\{dtq^k \, | \, k\in \mathbb{Z}: \, dtq^k<-1\}.
\end{split}
\end{equation}
We can extend the definition of the measure $m$ to a positive measure
for all parameters
$\beta$ satisfying \eqref{conditions}, using the fact that the
discrete weights $m\bigl(\{x\}\bigr)=m\bigl(\{x^{-1}\}\bigr)$ for $x\in D$
depend continuously on $\beta$, see \cite[(5.7) \& (5.8)]{KS2}.

Let $L^2_+(m)$ be the Hilbert space of
$L^2$-function $f$ with respect to the measure $m$
satisfying $f(x)=f(x^{-1})$  $m$-a.e.
The Askey-Wilson function transform $\mathcal{J}=\mathcal{J}(\beta;q)$
is now defined by
\begin{equation*}
\bigl(\mathcal{J}f\bigr)(\gamma)
=\frac{1}{K}\int f(x)\phi_\gamma(x)dm(x)
\end{equation*}
for compactly supported functions $f\in L_+^2(m)$,
with $K=K(\beta;q)$ the positive constant
\begin{equation}\label{Kon}
\begin{split}
K&=\frac{N^\tau}{\sqrt{\theta(adt,bdt,cdt,qt)}}\\
&=\frac{\bigl(qabc/d;q\bigr)_{\infty}}
{\bigl(q,ab,ac,bc,qa/d,qb/d,qc/d;q\bigr)_{\infty}}
\frac{1}{\sqrt{\theta(adt,bdt,cdt,qt)}},
\end{split}
\end{equation}
where $N$ is the Askey-Wilson integral, see \eqref{N}. 

We can now restate \cite[Theorem 1]{KS2} as follows.
\begin{thm}\label{mainprev}
The Askey-Wilson function transform $\mathcal{J}$
uniquely extends by continuity to an
isometric isomorphism
\[
\mathcal{J}: L^2_{+}(m)\rightarrow
L^2_{+}(m^\sigma).
\]
The inverse of $\mathcal{J}$ is given by
$\mathcal{J}^\sigma: L^2_{+}(m^\sigma)\rightarrow
L^2_{+}(m)$.
\end{thm}

Combined with \eqref{secondform}, we obtain the following
main result of this section.

\begin{thm}\label{master}
Suppose that the parameters $\beta=(\alpha,t)$ satisfy the conditions
\eqref{conditions}.
Let $s\in \mathcal{S}^\tau=\mathcal{S}^{\sigma\tau}$, then
\begin{equation}\label{masterformulas}
\begin{split}
\bigl(\mathcal{J}\bigl(E^\tau_sG^\tau\bigr)\bigr)(\gamma)&=
\frac{1}{\sqrt{\theta(adt,bdt,cdt,qt)}}
\frac{G^{\tau\sigma}(s_0^\tau)}{G^{\tau\sigma}(s)}E^{\sigma\tau}_s(\gamma)
G^{\sigma}(\gamma)^{-1}\Theta^{\sigma}(\gamma)^{-1},\\
\bigl(\mathcal{J}\bigl(E^{\tau}_sG^{-1}\Theta^{-1}\bigr)\bigr)(\gamma)&=
\sqrt{\theta(adt,bdt,cdt,qt)}\frac{G^{\sigma\tau\sigma}(s)}
{G^{\sigma\tau\sigma}(s_0^{\tau})}
E_s^{\sigma\tau}(\gamma)G^{\sigma\tau}(\gamma),
\end{split}
\end{equation}
as identities in $L^2_{+}(m^{\sigma})$.
\end{thm}
\begin{proof}
Observe that the factor $\theta(adt,bdt,cdt,qt)$ appearing in the
formulas is invariant under the parameter involution $\sigma$.
In view of Theorem \ref{mainprev} it thus suffices to prove
the explicit evaluation formula for
$\mathcal{J}\bigl(E_s^\tau G^{-1}\Theta^{-1}\bigr)$
with $s\in\mathcal{S}^\tau$. We show that it is in fact
a reformulation of \eqref{secondform}.

Observe that $E_s^\tau G^{-1}\Theta^{-1}\in L_+^2(m)$, since it
is a compactly supported function (use here that $\Theta^{-1}$ 
vanishes on the discrete mass points $D_-$ of the measure $m$).
For the moment we assume that the parameters $\beta=(\alpha,t)$
satisfy the conditions \eqref{conditions} and that they are generic.
Recall that the conditions \eqref{conditions} on the parameters $\beta$
imply that $0<b,c<1$ and $d>q$. By Cauchy's Theorem we thus conclude that
\begin{equation*}
\begin{split}
\frac{1}{4\pi i}\int_{\mathcal{T}^\tau}f(x)\Delta^\tau(x)\frac{dx}{x}
=\frac{1}{4\pi i}\int_{\mathbb{T}}&f(x)\Delta^\tau(x)\frac{dx}{x}\\
&+\frac{1}{2}\sum_{x\in D_+}\bigl(f(x)+f(x^{-1})\bigr)\,
\underset{y=x}{\hbox{Res}}\left(\frac{\Delta^\tau(y)}{y}\right)
\end{split}
\end{equation*}
for analytic, inversion-invariant functions $f: \mathbb{C}^\times \rightarrow
\mathbb{C}$, where $D_+$ is given by \eqref{D}.
Since $\Theta^{-1}$ vanishes on the discrete mass points $D_-$ of the
measure $m$, and since
\begin{equation*}
\begin{split}
\phi_\gamma(x)\Theta(x)^{-1}G(x)^{-1}W(x)&=
\phi_\gamma(x)G^\tau(x)^{-1}\Delta^\tau(x)\\
&=G^{\sigma\tau}(\gamma)\phi_\gamma^{an}(x)\Delta^\tau(x),
\end{split}
\end{equation*}
we obtain for $s\in \mathcal{S}^\tau$,
\begin{equation*}
\begin{split}
\bigl(\mathcal{J}\bigl(E^{\tau}_sG^{-1}\Theta^{-1}\bigr)\bigr)(\gamma)
&=\frac{1}{K}\left\{
\frac{1}{4\pi i}\int_{\mathcal{T}^\tau}E^\tau_s(x)
\phi_\gamma^{an}(x)\Delta^\tau(x)\frac{dx}{x}\right\}G^{\sigma\tau}(\gamma)\\
&=\frac{N^\tau}{K}\bigl(\mathbb{F}^\tau\bigl(\phi_\gamma^{an}\bigr)\bigr)(s)
G^{\sigma\tau}(\gamma)\\
&=\frac{G^{\tau\sigma\tau}(s)N^\tau}{G^{\tau\sigma\tau}(s_0^\tau)K}
E_s^{\sigma\tau}(\gamma)G^{\sigma\tau}(\gamma),
\end{split}
\end{equation*}
where the last equality is due to \eqref{secondform}.
The desired identity now follows from \eqref{parameterinvariant}
in view of the explicit expression \eqref{Kon} of the constant $K$.
The generic conditions on the parameters
$\beta$ can be removed by continuity.
\end{proof}

Consider the sub-spaces $V_{cl}=V_{cl}(\beta;q)$ and
$V_{str}=V_{str}(\beta;q)$ of $L_+^2(m)$ defined by
\[
V_{cl}=\mathcal{A} G^{-1}\Theta^{-1},\qquad
V_{str}=\mathcal{A} G^\tau.
\]
The subscripts
``cl'' and ``str'' stand for ``classical'' and ``strange'' respectively.
This terminology is justified by the harmonic analytic interpretation
of the Askey-Wilson function transform, see \cite{KS1}.
Indeed, when we regard $m$ as the Plancherel measure for the quantum
$SU(1,1)$ group, then the functions $f\in V_{cl}$ are supported
on the ``classical part'' of the measure $m$, which is the part of the
measure representing the contributions 
of the principal unitary series representations and of the positive
discrete series representations to the Plancherel measure. 
This translates to the property
that all functions $f\in V_{cl}$ vanish on the discrete mass points
$D_-$ of the measure $m$. The functions $V_{str}$ on the other
hand also have support on the ``strange part'' $D_-$ of the measure $m$,
which is the part of the measure representing the contributions 
of the strange series representations to the Plancherel measure.

Obviously $V_{cl}\cap V_{str}=\{0\}$ since any $f\in V_{cl}$
is zero on $D_-$.
Note that $V_{cl}\oplus V_{str}\subset L_+^2(m)$ is 
{\it not} an orthogonal direct sum
decomposition, since for all $s,u\in \mathcal{S}^\tau$,
\begin{equation*}
\int \frac{E_u^{\tau}(x)}{G(x)\Theta(x)}
E_s^\tau(x)G^\tau(x)dm(x)
=\frac{N^\tau}{h^{\tau\sigma}(s)}\delta_{s,u}
\end{equation*}
by the orthogonality relations \eqref{orthogonality}
for the Askey-Wilson polynomials, cf. the proof of Theorem
\ref{master}. In the next section it is shown 
that $V_{cl}\oplus V_{str}$ is
a dense sub-space of $L^2_+(m)$ if we impose an extra condition on the
allowed parameter values for $\beta$ 
(see Proposition \ref{density}). For these parameters $\beta$, the 
explicit formulas \eqref{masterformulas} thus completely determine 
the Askey-Wilson function transform. In particular, $V_{str}$ then 
completely takes care of the ``strange part'' $D_-$ of the measure $m$.

The formulas \eqref{masterformulas} immediately lead
to the following result.

\begin{cor}
The restriction
of the Askey-Wilson function transform $\mathcal{J}$ 
to the sub-space $V_{cl}\subset L_+^2(m)$ defines a linear bijection
$\mathcal{J}|_{V_{cl}}: V_{cl}\rightarrow V_{str}^\sigma$.
The inverse of $\mathcal{J}|_{V_{cl}}$ is given by 
$\mathcal{J}^\sigma|_{V^\sigma_{str}}$.
\end{cor}

The transform $\mathcal{J}|_{V_{cl}}$ is closely related to the
polynomial Askey-Wilson transform $\mathbb{F}^\tau$ acting upon
the sub-space $\mathcal{A}\,G^{\tau}{}^{-1}$. In fact, for generic
$\beta$ satisfying the conditions \eqref{conditions}, we have
\begin{equation}\label{linking}
\bigl(\mathcal{J}\bigl(E_s^\tau G^{-1}\Theta^{-1}\bigr)\bigr)(v)=
\frac{\bigl(qabc/d;q\bigr)_{\infty}N^\tau}
{\bigl(bc,q/ad,qa/d,qb/d,qc/d;q\bigr)_{\infty}K}
\bigl(\mathbb{F}^\tau\bigl(E_v\,G^{\tau}{}^{-1}\bigr)\bigr)(s)
\end{equation}
for all $v\in \mathcal{S}$ and $s\in \mathcal{S}^\tau$
by the polynomial reduction
\eqref{redpol} of the Askey-Wilson function, compare with 
the proof of Theorem \ref{master}. Observe in particular that 
the explicit formulas \eqref{polmasterformulas} for 
the polynomial Askey-Wilson transform
$\mathbb{F}^\tau$ acting on $\mathcal{A}\,G^{\tau}{}^{-1}$
are direct consequence of \eqref{linking}
and the explicit formulas \eqref{masterformulas} 
for $\mathcal{J}|_{V_{cl}}$.
 
On the other hand, by the strong convergence of the 
Gaussian $G_\tau(x)$ as $|x|$ tends
to zero it is possible to rewrite $\mathcal{J}|_{V_{str}}$
as a completely discrete transform by shrinking the radius of
the integration circle $\mathbb{T}$ in $m$ to zero while picking up
residues. 

These remarks show that the Askey-Wilson function transform 
$\mathcal{J}$ contains a continuous, 
polynomial (``compact'') type transform 
and a discrete, non-polynomial (``non-compact'') type transform
in a natural way, which are essentially each-others inverses. 
In the opposite direction, one may view the Askey-Wilson
function transform $\mathcal{J}$ as a {\it self-dual} 
transform obtained by glueing a continuous, compact type transform
and a discrete, non-compact type transform
together.

\begin{rem}
Cherednik's \cite{C} affine Hecke algebra approach 
extended to the present Askey-Wilson set-up
shows that there is a natural flexibility in the choice of the measure $m$.
More precisely, it turns out that for 
several different choices of measure $m$, 
the associated Fourier transform
$\mathcal{J}_m$ admits explicit
formulas which have a similar structure as the formulas
\eqref{masterformulas} for the Askey-Wilson function transform.  
This provides another explanation for
the similarities between Theorem \ref{polmaster} and
Theorem \ref{master}.
The ``proper choice'' of measure $m$ (and hence of 
Fourier transform $\mathcal{J}_m$) thus
depends on the applications which one has in mind.
For harmonic analysis on the quantum $SU(2)$ and quantum $SU(1,1)$
group, the ``proper choice'' is the polynomial Askey-Wilson transform
$\mathbb{F}$ and the Askey-Wilson function transform $\mathcal{J}$,
respectively (see \cite{K} and \cite{KS1}, respectively).
\end{rem}

We have used the inversion formula for
the Askey-Wilson function transform (see Theorem \ref{mainprev})
to explicitly compute the Askey-Wilson function transform
$\mathcal{J}$ acting upon functions $f\in V_{str}$, see Theorem \ref{master}.
On the other hand, the formulas \eqref{masterformulas} can be used to
reprove the inversion formula and Plancherel formula for the Askey-Wilson
function transform for those parameter values $\beta$
such that $V_{cl}\oplus V_{str}\subset L_+^2(m)$ is dense
(see Proposition \ref{density2} for the related density result).
It would therefore be of interest
to have an alternative proof of the formulas \eqref{masterformulas}
without referring to the main results of \cite{KS2}.
This is also of interest for the study of
multivariable generalizations of the present results.

An alternative proof for the key formulas \eqref{masterformulas}
is indeed possible using Cherednik's \cite{C}, \cite{C2}
theory of affine Hecke algebras, together with some elementary
elliptic function theory.
We do not go into details here, but only remark that these
techniques reduce the explicit computation of $\mathcal{J}$
on $V_{str}$ to the evaluation of $\bigl(\mathcal{J}(G^\tau)\bigr)(\tilde{a})$,
which is a $q$-analogue of the (one variable) Macdonald-Mehta
integral, cf. Macdonald \cite{Macdonald}. 
The integral $\bigl(\mathcal{J}(G^\tau)\bigr)(\tilde{a})$
can also be viewed as a ``non-polynomial'' analogue of the
(one variable) $q$-Macdonald-Mehta
integrals as discussed in Remark \ref{polMM}.
The evaluation of $\bigl(\mathcal{J}(G^\tau)\bigr)(\tilde{a})$
is equivalent to the following explicit identity.
\begin{thm}[One variable $q$-Macdonald-Mehta integral]\label{MMint}
For generic parameters $a,b,c,u\in \mathbb{C}^\times$ with
$q<|u|<1$ and $|a|,|b|,|c|<1$, we have
\begin{equation}\label{MacdonaldMehta}
\begin{split}
\frac{1}{4\pi i}&\int_{\mathbb{T}}\frac{\bigl(x^2,1/x^2;q\bigr)_{\infty}}
{\bigl(ax,a/x,bx,b/x,cx,c/x;q\bigr)_{\infty}\theta(ux,u/x)}\frac{dx}{x}\\
&+\sum_{k=1}^\infty\underset{x=uq^{-k}}{\textup{\hbox{Res}}}
\left(\frac{\bigl(x^2,1/x^2;q\bigr)_{\infty}}
{\bigl(ax,a/x,bx,b/x,cx,c/x;q\bigr)_{\infty}\theta(ux,u/x)x}\right)=\\
&\qquad\qquad\qquad\qquad\qquad\qquad\qquad\qquad=
\frac{1}{\bigl(q,ab,ac,bc;q\bigr)_{\infty}}
\frac{\theta(abcu)}{\theta(au,bu,cu)}.
\end{split}
\end{equation}
\end{thm}
\begin{proof}
For generic parameters $\beta$ satisfying \eqref{conditions}
and satisfying $|a|,|b|,|c|<1$ and $q<|dt|<1$,
we compute $\bigl(\mathcal{J}(G^\tau)\bigr)(\tilde{a})$ in two
different ways. The first way is by substituting the
definition of the Askey-Wilson function transform $\mathcal{J}$
and using the polynomial reduction \eqref{redpol}
for the Askey-Wilson function. This gives the
left hand side of
\eqref{MacdonaldMehta} with $u=dt$, multiplied by the constant
\[
\frac{\bigl(q,ab,ac;q\bigr)_{\infty}}
{\bigl(q/ad;q\bigr)_{\infty}}\sqrt{\theta(adt,bdt,cdt,qt)}.
\]
The second way to compute
$\bigl(\mathcal{J}(G^\tau)\bigr)(\tilde{a})$ is by using \eqref{masterformulas}
(take $s=s_0^\tau$ and $\gamma=\tilde{a}$ in the first formula
of \eqref{masterformulas}).
This gives the explicit infinite product evaluation
\[
\bigl(\mathcal{J}(G^\tau)\bigr)(\tilde{a})=
\frac{1}{\sqrt{\theta(adt,bdt,cdt,qt)}}\frac{\theta(qt,abcdt)}
{\bigl(bc,q/ad;q\bigr)_{\infty}}.
\]
Combining both expressions for
$\bigl(\mathcal{J}(G^\tau)\bigr)(\tilde{a})$
gives \eqref{MacdonaldMehta} with $u=dt$. The conditions on the parameters
can be removed by analytic continuation.

Due to its independent interest,
two alternative, direct proofs of \eqref{MacdonaldMehta}
are given in Appendix B.
\end{proof}

%%%%%%%%%%%%%%%%%%%%%%%%%%%%%%%%%%%%%%%%%%%%%%%%%%%%%%%%%%%%%%%%%
%%                                                             %%
%%          Density results                                    %%
%%                                                             %%
%%%%%%%%%%%%%%%%%%%%%%%%%%%%%%%%%%%%%%%%%%%%%%%%%%%%%%%%%%%%%%%%%

\section{Density results}

We next address the question whether $V_{cl}\oplus V_{str}$ is a dense
sub-space
of $L^2_+(m)$, i.e. whether the explicit formulas
\eqref{masterformulas} completely determine the Askey-Wilson function
transform. Surprisingly, the solution to this problem can be derived
from density results related to the polynomial Askey-Wilson transform
$\mathbb{F}$. We therefore first discuss the $L^2$-theory of
$\mathbb{F}$ and of its inverse $\mathbb{I}$.

For our purposes it suffices to restrict attention to
parameters $\alpha=(a,b,c,d)$ satisfying the conditions
\begin{enumerate}
\item[{\bf (i)}] $a,b,c,d\in\mathbb{R}$,
\item[{\bf (ii)}] $d<0<a$ and $abcd>0$,
\item[{\bf (iii)}] $ab,ac,ad,bc,bd,cd<1$.
\end{enumerate}
Under these conditions, at most two of the four parameters have
moduli $\geq 1$. If two parameters have moduli
$\geq 1$, then they have opposite sign.
For future reference, we write $V_{pol}$ for the set of four-tuples
$\alpha=(a,b,c,d)$ satisfying the conditions {\bf (i)}, {\bf (ii)} and
{\bf (iii)}.

For generic $\alpha\in V_{pol}$ and
for sufficiently regular, inversion-invariant functions $f$
(e.g. for $f\in \mathcal{A}$ or $f\in\mathcal{A}\,G^{-1}$),
we can rewrite the Fourier transform
$\bigl(\mathbb{F}f\bigr)(s)$\,
($s\in\mathcal{S}$) by Cauchy's Theorem as
\begin{equation}\label{F1}
 \bigl(\mathbb{F}f\bigr)(s)=\int f(x)E_s(x)d\nu(x),
\end{equation}
with $\nu=\nu(\cdot\,;\alpha;q)$ the positive measure
\begin{equation}\label{F2}
 \int f(x)d\nu(x)=\frac{1}{4\pi
i N}\int_{\mathbb{T}}f(x)\Delta(x)\frac{dx}{x}
+\frac{1}{2N}\sum_{x\in F}\bigl(f(x)+f(x^{-1})\bigr)\,
\underset{y=x}{\hbox{Res}}\left(\frac{\Delta(y)}{y}\right),
\end{equation}
with $F=F(\alpha;q)$ the finite, discrete set
\[ F=\{ eq^k \, | \, e\in \{a,b,c,d\},\,\, k\in \mathbb{Z}_+\,
\hbox{ such that } |eq^k|>1 \}
\]
and with $N$ the Askey-Wilson integral \eqref{N}.
By continuity in the parameters $\alpha$, both \eqref{F1} and
\eqref{F2} may be extended to all parameters $\alpha\in
V_{pol}$.
We use the notation $L^2_{+}\bigl(\nu\bigr)$ for the $L^2$-functions
$f: \mathbb{C}^\times\rightarrow \mathbb{C}$
with respect to the measure $\nu$ satisfying $f(x)=f(x^{-1})$
$\nu$-a.e. Sometimes it is convenient to think of the measure $\nu$ as a
positive measure $\widehat{\nu}$ on the real line supported on
\[ [-2,2]\cup \widehat{F},\qquad \widehat{F}=\{ x+x^{-1} \,\, | \,\, x\in F \},
\]
by the change of variable $y=x+x^{-1}$.
The measure $\widehat{\nu}$ is then given by
\[ \int f(y)d\widehat{\nu}(y)=\frac{1}{4\pi
N}\int_{-2}^2f(y)\frac{\widehat{\Delta}(y)}{(1-y^2/4)^{\frac{1}{2}}}dy
+\sum_{y\in\widehat{F}}f(y)\widehat{\nu}(\{y\}),
\]
with the weight function $\widehat{\Delta}$ 
satisfying $\widehat{\Delta}(x+x^{-1})=\Delta(x)$ and
with the discrete weights
\[ \widehat{\nu}(\{x+x^{-1}\}):=\nu(\{x\})+\nu(\{x^{-1}\})=2\nu(\{x\}),\qquad
x\in F.
\]
The Hilbert space
$L^2_{+}(\nu)$ can then be identified with the Hilbert space
$L^2(\widehat{\nu})$ of $L^2$-functions with respect
to the measure $\widehat{\nu}$ (cf. \cite{AW}). Observe that
the space $\mathbb{C}[x]$ of polynomials with complex coefficients
is dense in $L^2(\widehat{\nu})$ since $\widehat{\nu}$ is compactly supported.
Equivalently, $\mathcal{A}\subset L^2_+(\nu)$ is a dense sub-space.

For $\alpha\in V_{pol}$ we define the discrete measure 
$\mu=\mu(\cdot\,;\alpha;q)$ supported on $\mathcal{S}$ by
\[ \int g(x)d\mu(x)= \sum_{s\in \mathcal{S}}
g(s)\widetilde{h}(s).
\]
Observe that $\mu$ is a positive measure since the inverse quadratic norm
$\widetilde{h}(s)$\, (see \eqref{hnorm}) of the Askey-Wilson
polynomial $E_s$ is strictly positive for
$\alpha\in V_{pol}$.
Let $L^2(\mu)$ be the corresponding $L^2$-space.
By continuity, Theorem \ref{AWresult} implies the following
result.
\begin{cor}\label{L2result}
Let $\alpha\in V_{pol}$. The polynomial Askey-Wilson
transform $\mathbb{F}$ and its inverse $\mathbb{I}$ extend
by continuity to isometric isomorphisms
\begin{equation*}
\begin{split}
\mathbb{F}:\, &L^2_{+}\bigl(\nu\bigr)\rightarrow L^2(\mu),\\
\mathbb{I}:\, &L^2(\mu)\rightarrow L^2_{+}\bigl(\nu\bigr).
\end{split}
\end{equation*}
Furthermore, $\mathbb{I}: L^2(\mu)\rightarrow
L^2_{+}\bigl(\nu\bigr)$ is the inverse of
$\mathbb{F}: L^2_{+}\bigl(\nu\bigr)\rightarrow L^2(\mu)$.
\end{cor}

Combined with Theorem \ref{polmaster} we obtain  
the following lemma.

\begin{lem}\label{densitylemma}
Let $\alpha\in V_{pol}$. The orthocomplement of $\mathcal{A}\,G^{-1}$ in
$L_+^2(\nu)$ \textup{(}respectively of
$\bigl(\mathcal{A}\,G^{\sigma\tau}\bigr)|_{\mathcal{S}}$ 
in $L^2(\mu)$\textup{)} is finite dimensional. 
In both cases, the dimension of the orthocomplement is
\[\#\{ n\in \mathbb{Z}_+ \,\, | \,\,\, |dq^n|>1 \}.
\]
In particular, $\mathcal{A}\,G^{-1}\subset L_+^2(\nu)$
\textup{(}respectively $\bigl(\mathcal{A}\,G^{\sigma\tau}\bigr)|_{\mathcal{S}}
\subset L^2(\mu)$\textup{)} is dense if and only if $|d|\leq 1$.
\end{lem}
\begin{proof}
Let $\alpha=(a,b,c,d)\in V_{pol}$. Let $\widehat{G}$ be the
meromorphic function
satisfying $\widehat{G}(x+x^{-1})=G(x)$.
Let $\widehat{F}_0$ be the intersection of the polar divisor of
$\widehat{G}$ with the set $\widehat{F}$ of discrete mass points of
$\widehat{\nu}$. It is easy to verify that
\[
 \widehat{F}_0=
\{ dq^n+d^{-1}q^{-n} \,\, | \,\, n\in \mathbb{Z}_+:\,\, |dq^n|>1 \}.
\]
The closure of
$\mathbb{C}[x]\,\widehat{G}^{-1}$ in $L^2(\widehat{\nu})$ 
is exactly the sub-space of functions
$f\in L^2(\widehat{\nu})$ which vanish on $\widehat{F}_0$,
since $\widehat{\nu}$ is compactly
supported. Hence the orthocomplement of $\mathbb{C}[x]\,G^{-1}$ 
in $L^2(\widehat{\nu})$ 
is a $\#\widehat{F}_0$-dimensional sub-space of $L^2(\widehat{\nu})$. 
Equivalently, the orthocomplement of $\mathcal{A}\,G^{-1}$
in $L_+^2(\nu)$ is a $\#\widehat{F}_0$-dimensional sub-space
of $L_+^2(\nu)$. 

The identities in Theorem \ref{polmaster}, which were proven for
generic parameters $\alpha$, are valid for all $\alpha\in V_{pol}$ 
since they are regular at $\alpha\in V_{pol}$. 
By Theorem \ref{polmaster} and Corollary
\ref{L2result} the previous results on $\mathcal{A}\,G^{-1}\subset
L_+^2(\nu)$ thus imply that the orthocomplement of
$\bigl(\mathcal{A}\,G^{\sigma\tau}\bigr)|_{\mathcal{S}}$ in
$L^2(\mu)$ is a $\#\widehat{F}_0$-dimensional sub-space of $L^2(\mu)$.
\end{proof}

\begin{rem}
Suppose that $\alpha\in V_{pol}$ with $|d|>1$. It follows from 
Theorem \ref{L2result} and from the proof of Lemma \ref{densitylemma} 
that the functions $f_n\in L^2(\mu)$\, 
($n\in\mathbb{Z}_+$: $|dq^n|>1$) defined by
\[ f_{n}(s)=E_s(dq^n)=
{}_4\phi_3\left(\begin{matrix} \tilde{a}s, \tilde{a}/s, adq^n, q^{-n}a/d\\
ab, ac, ad\end{matrix}\,; q,q\right),\qquad s\in\mathcal{S}
\]
form an orthogonal basis for the orthocomplement of 
$\bigl(\mathcal{A}G^{\sigma\tau}\bigr)|_{\mathcal{S}}$ in $L^2(\mu)$.
The quadratic norm of $f_n$\, ($n\in\mathbb{Z}_+$: $|dq^n|>1$)
in $L^2(\mu)$ is given by
\[ \int |f_n(s)|^2d\mu(s)=
\frac{N}{\underset{y=dq^n}{\hbox{Res}}\left(\frac{\Delta(y)}{y}\right)}
\]
for the generic parameters $\alpha\in V_{pol}$ such that the pole 
of $\Delta(y)$ at $y=dq^n$ is simple. 
\end{rem}

\begin{cor}\label{density}
Let $u,v\in \mathbb{R}^\times$ with $|u|\leq 1$ and $v>0$.
Let $\rho=\rho_{u,v}$ be the positive discrete measure given by
\begin{equation}\label{rho}
 \int f(x)d\rho(x)=\sum_{k=0}^{\infty}
f(uq^k+u^{-1}q^{-k})v^{k}q^{k(k-1)}
\end{equation}
and let $L^2(\rho)$ be the associated $L^2$-space. 
Then $\mathbb{C}[x]$ is dense
in $L^2(\rho_{u,v})$ if and only if $v\leq 1$.
\end{cor}
\begin{proof}
Let $\alpha=(a,b,c,d)\in V_{pol}$. Observe that the condition 
$\alpha\in V_{pol}$ implies $\tilde{a}<q^{-\frac{1}{2}}$, so 
$s_n+s_n^{-1}=s_m+s_m^{-1}$ for 
$m,n\in\mathbb{Z}_+$ iff $n=m$, where (recall) $s_n=\tilde{a}q^n$. 
Hence we can define a positive, discrete measure $\widehat{\rho}$ supported on
$s_n+s_n^{-1}$\, ($n\in\mathbb{Z}_+$),
with weights
\[ \widehat{\rho}(\{s_n+s_n^{-1}\})=|G^{\sigma\tau}(s_n)|^2\mu(\{s_n\}),\qquad
n\in\mathbb{Z}_+.
\]
The corresponding $L^2$-space $L^2(\widehat{\rho})$ is
isomorphic to $L^2(\mu)$ via the surjective isometric isomorphism
$T: L^2(\mu)\rightarrow L^2(\widehat{\rho})$ defined by
\[\bigl(Tf\bigr)(s_n+s_n^{-1})=G^{\sigma\tau}(s_n)^{-1}f(s_n),\qquad
n\in\mathbb{Z}_+.
\]
Furthermore, the image of 
$\bigl(\mathcal{A}\,G^{\sigma\tau}\bigr)|_{\mathcal{S}}$  under
$T$ is exactly the space $\mathbb{C}[x]\subset L^2(\widehat{\rho})$
of polynomials with complex coefficients.
It thus follows from the previous lemma
that $\mathbb{C}[x]$ is dense in $L^2(\widehat{\rho})$ iff $|d|\leq 1$.

A direct computation using \eqref{hnorm} shows that
\[ \widehat{\rho}(\{s_n+s_n^{-1}\})=
\gamma(s_n+s_n^{-1})\rho_{\tilde{a},d^2}(\{s_n+s_n^{-1}\}), 
\qquad n\in \mathbb{Z}_+
\]
with $\rho_{u,v}$ the positive, discrete measure defined by \eqref{rho} and 
with 
\[\gamma: \hbox{supp}(\widehat{\rho})=\hbox{supp}(\rho_{\tilde{a},d^2})
\rightarrow \mathbb{R}_{>0}
\]
 a strictly positive,
bounded function with bounded inverse.  We conclude that
$\mathbb{C}[x]$ is dense in $L^2(\rho_{\tilde{a},d^2})$ iff $|d|\leq 1$.

Choose now $v>0$ and $0<u\leq 1$ arbitrarily.
Then there exist parameters $\alpha=(a,b,c,d)\in V_{pol}$ 
such that
\[ d^2=v,\qquad \tilde{a}=u.
\]
This proves the corollary in case $0<u\leq 1$. The corollary for $-1\leq u<0$
follows from the result for $0<u\leq 1$, using the surjective
isometric isomorphism $S: L^2(\rho_{u,v})\rightarrow L^2(\rho_{-u,v})$
defined by $(Sf)(x)=f(-x)$.
\end{proof}

\begin{rem}
Borichev and Sodin \cite[Theorem A]{BS} formulated 
criteria  
for the density of the space of polynomials 
$\mathbb{C}[x]$ in a $L^p$-space ($p\geq 1$)
when the associated measure is supported on the zero set of a 
Hamburger class function $B$, see \cite{BS} for more details.
The measure $\rho_{u,v}$ is of this particular form, with the 
associated Hamburger class function $B$ given by
 \[ B(z)=\prod_{\lambda\in \hbox{supp}(\rho_{u,v})}
\left(1-\frac{z}{\lambda}\right).
\]
It is a nice exercise to re-prove Corollary \ref{density}  
using these general criteria of Borichev and Sodin.
\end{rem}
\begin{rem}
By a result of M. Riesz (see e.g. \cite[Lemma A]{BD}), 
it is easy to verify that the measure 
$\rho_{u,v}$ for $u,v\in\mathbb{R}^\times$ with 
$|u|\leq 1$ and $v>0$ corresponds to a determinate moment problem
if and only if $0<v\leq q^2$. 
In particular, $\rho_{u,v}$ is a $N$-extremal measure 
(or, in the terminology of \cite{BS}, a canonical measure) if and only
if $q^2<v\leq 1$. If $\rho_{u,v}$ is determinate (i.e. if $v\leq q^2$),
then $\rho_{u,v}$ has a finite index of determinacy, 
which can be computed explicitly (see Berg and Duran \cite{BD}
for a detailed study of measures with finite index of determinacy). 
On the other hand, if $\rho_{u,v}$
is indeterminate and the polynomials are not dense in
$L^2(\rho_{u,v})$ (i.e. if $v>1$), then the dimension of the
orthocomplement of $\mathbb{C}[x]$ in $L^2(\rho_{u,v})$ equals
$\#\{ k\in \mathbb{Z}_+ \, | \,\, |vq^{2k}|>1 \}$ 
by e.g. \cite[Proposition A1.4]{BS}. This observation nicely relates to 
Lemma \ref{densitylemma}.
\end{rem}

We now use
Corollary \ref{density} to derive the following density result for the
linear sub-space $V_{cl}\oplus V_{str}\subset L_+^2(m)$.

\begin{prop}\label{density2}
Let the parameters $\beta=(\alpha,t)$ satisfy the conditions
\eqref{conditions}. Let $k\in\mathbb{Z}$ be the unique integer
such that $1<|dtq^k|\leq q^{-1}$.

The sub-space
\[ V_{cl}\oplus V_{str} \subset L_{+}^2(m)
\]
is dense if and only if
$|\tilde{a}tq^k|\geq 1$. Furthermore, the \textup{(}non-empty\textup{)}
set of parameters
$\beta$ satisfying \eqref{conditions} and satisfying the condition
 $|\tilde{a}tq^k|\geq 1$ is invariant under the duality involution $\sigma$.
\end{prop}
\begin{proof}
We first prove the last part of the proposition. Suppose that 
the parameters $\beta$
satisfy \eqref{conditions} and $|\tilde{a}tq^k|\geq 1$, where $k\in
\mathbb{Z}$ is the unique integer such that $1<|dtq^k|\leq q^{-1}$.
The dual parameters $\widetilde{\beta}=\beta_\sigma$ then 
satisfy the conditions \eqref{conditions} in view of \cite[Lemma
1]{KS2}. It remains to verify the inequality $|a\tilde{t}q^r|\geq 1$,
where $r\in\mathbb{Z}$ is the unique integer such that
$1<|\tilde{d}\tilde{t}q^{r}|\leq q^{-1}$. 
By the definition of dual parameters the condition
$1<|\tilde{d}\tilde{t}q^{r}|\leq q^{-1}$ is equivalent to the condition
$q\leq |\tilde{a}tq^{-r}|<1$. In particular, $-r>k$. But then
\[ |a\tilde{t}q^r|\geq
q^{-1}|a\tilde{t}q^{-k}|=q^{-1}|d^{-1}t^{-1}q^{-k}|\geq 1,
\]
which is the desired inequality.

We now focus on the first part of the statement. We fix
parameters $\beta$ satisfying the conditions \eqref{conditions}.
Via the change of variable $y=x+x^{-1}$ we can rewrite the measure
$m$ as a positive measure $\widehat{m}$ on $\mathbb{R}$ supported
on $\mathcal{D}_{cl}\cup\mathcal{D}_{str}$, where
\begin{equation*}
\begin{split}
\mathcal{D}_{cl}&=[-2,2]\cup \{aq^n+a^{-1}q^{-n} \,\, | \,\,
n\in\mathbb{Z}_+ :\,\, aq^n>1 \},\\
\mathcal{D}_{str}&=\{ uq^{n}+u^{-1}q^{-n} \,\, | \,\,
n\in\mathbb{Z}_+\},
\end{split}
\end{equation*}
with $u=d^{-1}t^{-1}q^{-k}$.
Here $k\in\mathbb{Z}$ is the unique integer 
such that $1<|dtq^k|\leq q^{-1}$. Under the change
of variable $y=x+x^{-1}$, the Hilbert space $L_{+}^2(m)$ is isomorphic to
the Hilbert space $L^2(\widehat{m})$ of $L^2$-functions
with respect to the measure $\widehat{m}$ (compare with the 
identification of $L^2_+(\nu)$ and
$L^2(\widehat{\nu})$ as discussed at the beginning of this section). 

Consider $\widehat{V}_{cl}=\mathbb{C}[x]\,\varphi$ and  
$\widehat{V}_{str}=\mathbb{C}[x]\,\widehat{G}^\tau$
as linear sub-spaces of $L^2(\widehat{m})$, where
$\varphi$ and $\widehat{G}^\tau$ are the meromorphic functions satisfying
\[ \varphi(x+x^{-1})=\Theta(x)^{-1}G(x)^{-1},\qquad
\widehat{G}^\tau(x+x^{-1})=G^\tau(x).
\]
Then $V_{cl}\oplus V_{str}\subset L_+^2(m)$ is dense iff
$\widehat{V}_{cl}\oplus \widehat{V}_{str}\subset L^2(\widehat{m})$
is dense.

Let $\widehat{m}_{cl}=\widehat{m}|_{\mathcal{D}_{cl}}$ 
(respectively $\widehat{m}=\widehat{m}|_{\mathcal{D}_{str}}$)
be the restriction of the measure $\widehat{m}$
to $\mathcal{D}_{cl}$ (respectively $\mathcal{D}_{str}$), and denote
$L^2(\widehat{m}_{cl})$ (respectively $L^2(\widehat{m}_{str})$)
for the associated $L^2$-space. We define surjective,
continuous linear mappings
\[\pi_{cl}: L^2(\widehat{m})\rightarrow L^2(\widehat{m}_{cl}),\qquad
\pi_{str}: L^2(\widehat{m})\rightarrow L^2(\widehat{m}_{str})
\]
by $\pi_{cl}(f)=f|_{\mathcal{D}_{cl}}$ and
$\pi_{str}(f)=f|_{\mathcal{D}_{str}}$. 

Observe that $\pi_{cl}(\varphi)$ is non-zero $\widehat{m}_{cl}$-a.e.
due to the conditions \eqref{conditions} on the parameters $\beta$.
Since the measure $\widehat{m}_{cl}$ is compactly 
supported, we conclude that the sub-space $\pi_{cl}(\widehat{V}_{cl})
\subset L^2(\widehat{m}_{cl})$ is dense,
compare with the proof of Lemma \ref{densitylemma}. 

Let $H_{cl}\subset L^2(\widehat{m})$ be the closed sub-space of
functions $f\in L^2(\widehat{m})$ with support contained in
$\mathcal{D}_{cl}$. Then $\widehat{V}_{cl}\subset H_{cl}$ since $\varphi$
vanishes on $\mathcal{D}_{str}$, and $\pi_{cl}|_{H_{cl}}:
H_{cl}\rightarrow L^2(\widehat{m}_{cl})$ is a surjective isometric
isomorphism. It follows that $\widehat{V}_{cl}\subset H_{cl}$ is dense.

Since $\widehat{V}_{cl}\subset H_{cl}$ is dense we have that
$\widehat{V}_{cl}\oplus \widehat{V}_{str}\subset 
L^2(\widehat{m})$ is dense iff 
$\pi_{str}(\widehat{V}_{str})\subset L^2(\widehat{m}_{str})$ is dense.
It thus suffices to prove that $\pi_{str}(\widehat{V}_{str})\subset
L^2(\widehat{m}_{str})$ is dense iff the parameters $\beta$
satisfy the extra condition $|\tilde{a}tq^k|\geq 1$.
Observe first that $-1<u=d^{-1}t^{-1}q^{-k}<0$ by the conditions
\eqref{conditions} on the parameters $\beta$ and by the definition of
the integer $k$. 
Furthermore, for any discrete mass point
$y_n=uq^n+u^{-1}q^{-n}\in\mathcal{D}_{str}$\, ($n\in\mathbb{Z}_+$),
we have
\[ |\widehat{G}^\tau(y_n)|^2\widehat{m}_{str}(\{y_n\})=
\gamma(y_n)\rho_{u,v}(\{y_n\}),\qquad v:=\tilde{a}^{-2}t^{-2}q^{-2k},
\]
with $\rho=\rho_{u,v}$ the measure \eqref{rho}
and with $\gamma: \mathcal{D}_{str}\rightarrow \mathbb{R}_{>0}$ a 
bounded function with bounded inverse 
(we have use here the explicit expression for
the weights $\widehat{m}_{str}(\{y_n\})$, see \cite[(5.8)]{KS2}).
It follows that $\pi_{str}(\widehat{V}_{str})\subset 
L^2(\widehat{m}_{str})$ is dense iff 
$\mathbb{C}[x]\subset L^2(\rho_{u,v})$ is dense, cf. the proof of
Corollary \ref{density}.
The desired density result is now a direct consequence of 
Corollary \ref{density}.
\end{proof}

%%%%%%%%%%%%%%%%%%%%%%%%%%%%%%%%%%%%%%%%%%%%%%%%%%%%%%%%%%%%%%%%%%
%%                                                              %%
%%   Appendix A: Proof of the expansion formula                 %%
%%                                                              %%
%%%%%%%%%%%%%%%%%%%%%%%%%%%%%%%%%%%%%%%%%%%%%%%%%%%%%%%%%%%%%%%%%%

\section{Appendix A: Proof of Proposition \ref{expansion2}}

In this appendix we give a proof of the formulas \eqref{secondform},
which can be rewritten as
\begin{equation}\label{newform}
\begin{split}
\bigl(\mathbb{F}^\tau\bigl(\phi_\gamma^{an}\bigr)\bigr)(s_m^\tau)=
&\frac{\bigl(bc;q\bigr)_m}{\bigl(qa/d;q\bigr)_m}\left(\frac{-qa}{d}\right)^m
q^{m(m-1)/2}\\
&\quad\qquad\times{}_4\phi_3\left(\begin{matrix} q^{-m}, q^mabc/d,
\tilde{a}\gamma, \tilde{a}/\gamma\\
ab,ac,bc\end{matrix}\,; q,q\right)
\end{split}
\end{equation}
for $m\in\mathbb{Z}_+$.

Throughout the proof of \eqref{newform} we fix
$m\in\mathbb{Z}_+$, and we set $s=s_m^\tau$.
We substitute the expression \eqref{43presentation} for
the Askey-Wilson function $\phi_\gamma(x)$
in the integral
\begin{equation*}
\begin{split}
 \bigl(\mathbb{F}^\tau\bigl(\phi_\gamma^{an}\bigr)\bigr)(s)&=
\frac{1}{4\pi i N^\tau}\int_{\mathcal{T}^\tau}
\phi_\gamma^{an}(x)E_{s}^\tau(x)\Delta^\tau(x)\frac{dx}{x}\\
&=\frac{1}{4\pi i N^\tau G^{\sigma\tau}(\gamma)}\int_{\mathcal{T}^\tau}
\phi_\gamma(x)E_{s}^\tau(x)
\frac{\Delta^\tau(x)}{G^\tau(x)}\frac{dx}{x}
\end{split}
\end{equation*}
and we use that
\[
E_{s}^\tau(x)=
\frac{\bigl(bc,qb/d;q\bigr)_m}{\bigl(ac,qa/d;q\bigr)_m}
\left(\frac{a}{b}\right)^m{}_4\phi_3
\left(\begin{matrix} q^{-m}, q^mabc/d, bx, b/x\\
 ab, bc, qb/d \end{matrix}\,; q,q\right)
\]
by Sear's transformation formula \cite[(2.10.4)]{GR} with $a,b,c,d,e$ and $f$
in \cite[(2.10.4)]{GR} taken to be $q^mabc/d$, $ax$, $a/x$, $ab$, $ac$ and
$qa/d$, respectively.
We arrive at
\begin{equation}\label{hulp1}
\bigl(\mathbb{F}^\tau\bigl(\phi_\gamma^{an}\bigr)\bigr)(s)
=\frac{\bigl(bc,qb/d;q\bigr)_m}{\bigl(ac,qa/d;q\bigr)_m}
\left(\frac{a}{b}\right)^m\left\{I_1(\gamma)+I_2(\gamma)\right\},
\end{equation}
with $I_1(\gamma)$ given by
\begin{equation*}
\begin{split}
I_1(\gamma)=&\frac{\bigl(qabc/d,q\gamma/\tilde{d},
q/\tilde{d}\gamma;q\bigr)_{\infty}}
{\bigl(bc,qa/d,qb/d,qc/d,q/ad;q\bigr)_{\infty}}\\
&\qquad\times
\sum_{n=0}^m\sum_{k=0}^\infty
\frac{\bigl(\tilde{a}\gamma,\tilde{a}/\gamma;q\bigr)_k
\bigl(q^{-m},q^mabc/d;q\bigr)_n}
{\bigl(q,ab,ac,ad;q\bigr)_k\bigl(q,ab,bc,qb/d;q\bigr)_n}q^{k+n}\\
&\qquad\qquad\qquad\times
\frac{1}{4\pi iN^\tau}\int_{\mathcal{T}^\tau}
\frac{\bigl(x^2,1/x^2;q\bigr)_{\infty}}
{\bigl(q^kax,q^ka/x,q^nbx,q^nb/x,cx,c/x;q\bigr)_{\infty}}\frac{dx}{x},
\end{split}
\end{equation*}
and with $I_2(\gamma)$ given by
\begin{equation*}
\begin{split}
I_2(\gamma)=&\frac{\bigl(qabc/d,\tilde{a}\gamma,
\tilde{a}/\gamma;q\bigr)_{\infty}}
{\bigl(ab,ac,bc,qa/d,ad/q;q\bigr)_{\infty}}\\
&\qquad\times
\sum_{n=0}^m\sum_{k=0}^\infty
\frac{\bigl(q\gamma/\tilde{d},q/\tilde{d}\gamma;q\bigr)_k
\bigl(q^{-m},q^mabc/d;q\bigr)_n}
{\bigl(q,qb/d,qc/d,q^2/ad;q\bigr)_k\bigl(q,ab,bc,qb/d;q\bigr)_n}q^{k+n}\\
&\qquad\qquad\times
\frac{1}{4\pi iN^\tau}\int_{\mathcal{T}^\tau}
\frac{\bigl(x^2,1/x^2;q\bigr)_{\infty}}
{\bigl(q^nbx,q^nb/x,cx,c/x,q^{k+1}x/d,q^{k+1}/dx;q\bigr)_{\infty}}\frac{dx}{x}.
\end{split}
\end{equation*}
Now the integrals in the expressions for $I_1(\gamma)$ and
$I_2(\gamma)$ can be evaluated as special case of
the evaluation of the Askey-Wilson integral,
see \eqref{N}.
The resulting sum over $k$ in both $I_1(\gamma)$ and $I_2(\gamma)$
can then be rewritten as a non-terminating ${}_3\phi_2$.
This leads to the identity
\begin{equation*}
\begin{split}
&I_1(\gamma)+I_2(\gamma)=\\
&=\sum_{n=0}^m
\frac{\bigl(q^{-m},
q^mabc/d;q\bigr)_n}{\bigl(q;q\bigr)_n\bigl(bc;q\bigr)_{\infty}}q^n
\left\{\frac{\bigl(q\gamma/\tilde{d},q/\tilde{d}\gamma;q\bigr)_{\infty}}
{\bigl(q/ad;q\bigr)_{\infty}\bigl(qb/d;q\bigr)_n}
{}_3\phi_2
\left(\begin{matrix} \tilde{a}\gamma,\tilde{a}/\gamma,q^nab\\
 ab,ad\end{matrix}\,; q,q\right)\right.\\
&\qquad\qquad\qquad\qquad\qquad
\left.+\frac{\bigl(\tilde{a}\gamma,\tilde{a}/\gamma;q\bigr)_{\infty}}
{\bigl(ad/q;q\bigr)_{\infty}\bigl(ab;q\bigr)_n}
{}_3\phi_2
\left(\begin{matrix} q\gamma/\tilde{d},q/\tilde{d}\gamma,q^{n+1}b/d\\
 q^2/ad,qb/d\end{matrix}\,; q,q\right)\right\}.
\end{split}
\end{equation*}
Applying the three-term transformation formula
\cite[(3.3.1)]{GR} for ${}_3\phi_2$'s with parameters $a,b,c,d$ and $e$ in
\cite[(3.3.1)]{GR} specialized to
$q^{-n}$, $\tilde{a}\gamma$,
$\tilde{a}/\gamma$, $ab$ and
$bc(=q\tilde{a}/\tilde{d})$ respectively, shows that
\begin{equation}\label{hulp2}
I_1(\gamma)+I_2(\gamma)=
\sum_{n=0}^m\frac{\bigl(q^{-m},q^mabc/d;q\bigr)_n}
{\bigl(q,qb/d;q\bigr)_n}q^n
{}_3\phi_2
\left(\begin{matrix} q^{-n},\tilde{a}\gamma,\tilde{a}/\gamma\\
 ab,bc\end{matrix}\,; q,q^{n+1}b/d\right).
\end{equation}
Formula \eqref{newform} now immediately follows from \eqref{hulp1},
\eqref{hulp2} and the following lemma.
\begin{lem}\label{hulp3}
The following identity is valid:
\begin{equation*}
\begin{split}
\sum_{n=0}^m&\frac{\bigl(q^{-m},q^mabc/d;q\bigr)_n}
{\bigl(q,qb/d;q\bigr)_n}q^n
{}_3\phi_2
\left(\begin{matrix} q^{-n},\tilde{a}\gamma,\tilde{a}/\gamma\\
 ab,bc\end{matrix}\,; q,q^{n+1}b/d\right)=\\
&=\left(\frac{-qb}{d}\right)^mq^{m(m-1)/2}\frac{\bigl(ac;q\bigr)_m}
{\bigl(qb/d;q\bigr)_m}
{}_4\phi_3\left(\begin{matrix} q^{-m},q^mabc/d,
\tilde{a}\gamma,\tilde{a}/\gamma\\
ab,ac,bc\end{matrix}\,; q,q\right).
\end{split}
\end{equation*}
\end{lem}
\begin{proof}
Denote the left hand side of the desired identity by $f_m(\gamma)$.
It is clear that $f_m(\gamma)$ is a polynomial of degree $m$
in $\gamma+\gamma^{-1}$. In the expansion
\[
f_m(\gamma)=
\sum_{k=0}^m\alpha_k\bigl(\tilde{a}\gamma,\tilde{a}/\gamma;q\bigr)_k,
\]
the coefficients $\alpha_k$ are explicitly given by
\[
\alpha_k=
\sum_{n=k}^m\frac{\bigl(q^{-m},q^mabc/d;q\bigr)_n\bigl(q^{-n};q\bigr)_k}
{\bigl(q,qb/d;q\bigr)_n\bigl(q,ab,bc;q\bigr)_k}q^n\left(\frac{q^{n+1}b}{d}
\right)^k.
\]
Changing the summation variable to $r=n-k$ and simplifying the sum
yields
\begin{equation*}
\begin{split}
\alpha_k=\frac{\bigl(q^{-m},q^mabc/d;q\bigr)_k}
{\bigl(q,q,ab,bc,qb/d;q\bigr)_k}&\Bigl(\frac{q^{k+2}b}{d}\Bigr)^k\\
\times&\sum_{r=0}^{m-k}
\frac{\bigl(q^{k-m},q^{m+k}abc/d;q\bigr)_r}
{\bigl(q^{k+1},q^{k+1}b/d;q\bigr)_r}\bigl(q^{-r-k};q\bigr)_kq^{(k+1)r}.
\end{split}
\end{equation*}
Now applying \cite[(1.2.37)]{GR} to the $q$-shifted factorial
$\bigl(q^{-r-k};q\bigr)_k$ appearing in the right hand side of the last
formula for $\alpha_k$, leads to
\[
\alpha_k=\frac{\bigl(q^{-m},q^mabc/d;q\bigr)_k}
{\bigl(q,ab,bc,qb/d;q\bigr)_k}q^{k(k-1)/2}\left(\frac{-q^2b}{d}\right)^k
{}_2\phi_1\left(\begin{matrix} q^{k-m},q^{m+k}abc/d\\
q^{k+1}b/d\end{matrix}\,; q,q\right).
\]
The terminating ${}_2\phi_1$ is summable by the $q$-Vandermonde
formula \cite[(1.5.3)]{GR}. Simplification of the resulting expression
then shows that
\[
\alpha_k=\left(\frac{-qb}{d}\right)^mq^{m(m-1)/2}
\frac{\bigl(ac;q\bigr)_m}{\bigl(qb/d;q\bigr)_m}
\frac{\bigl(q^{-m},q^mabc/d;q\bigr)_k}
{\bigl(q,ab,ac,bc;q\bigr)_k}q^k,
\]
as desired.
\end{proof}

%%%%%%%%%%%%%%%%%%%%%%%%%%%%%%%%%%%%%%%%%%%%%%%%%%%%%%%%%%%%%%%%%%%%
%%                                                                %%
%%   Appendix B: Evaluation of the one variable                   %%
%%   $q$-Macdonald-Mehta integral                                 %%
%%                                                                %%
%%%%%%%%%%%%%%%%%%%%%%%%%%%%%%%%%%%%%%%%%%%%%%%%%%%%%%%%%%%%%%%%%%%%

\section{Appendix B: Evaluations of the one variable
$q$-Macdonald-Mehta integral}
In this appendix we give two alternative proofs of the 
$q$-analogue of the (one variable) Macdonald-Mehta integral, 
see Theorem \ref{MMint}. The first proof is
based on Nassrallah's and Rahman's \cite{NR}
integral representation of the very-well-poised ${}_8\phi_7$ series.
The second proof uses the fact that the $q$-analogue of the
Macdonald-Mehta integral in one variable can be rewritten in
a completely discrete form using Cauchy's Theorem.
The evaluation then follows from
limit cases of the summation formulas of the very-well-poised
${}_6\phi_5$ series and of the very-well-poised ${}_6\psi_6$ series,
together with some elementary elliptic function theory (the second
proof is in the spirit of Askey's and Wilson's \cite{AW} original proof of
the evaluation of the Askey-Wilson integral).

\subsection{First direct proof of Theorem \ref{MMint}}
We fix generic parameters $a,b,c,u\in\mathbb{C}^\times$ satisfying
$|abcu|>q$ and $|a|,|b|,|c|,|u|<1$. We write
\[ L_1=\frac{1}{4\pi i}\int_{\mathbb{T}}
\frac{\bigl(x^2,1/x^{2};q\bigr)_{\infty}}
{\bigl(ax,a/x,bx,b/x,cx,c/x;q\bigr)_{\infty}\,\theta(ux,u/x)}\frac{dx}{x}
\]
and
\[ L_2=\sum_{k=1}^\infty\underset{x=uq^{-k}}{\hbox{Res}}
\left(\frac{\bigl(x^2,1/x^{2};q\bigr)_{\infty}}
{\bigl(ax,a/x,bx,b/x,cx,c/x;q\bigr)_{\infty}\,\theta(ux,u/x)x}\right),
\]
respectively. The aim is to evaluate $L_1+L_2$.
For $L_1$, we use the integral representation \cite[(6.3.8)]{GR}
of the very-well-poised ${}_8\phi_7$
series due to Nassrallah and Rahman \cite{NR}, with
parameters $a,b,c,d,f$ and $g$ in \cite[(6.3.8)]{GR} specialized to
$a,b,c,u,q/u$ and $0$, respectively. The ${}_8W_7$ then reduces to a
${}_3\phi_2$, and we obtain
\begin{equation}\label{stepB1}
L_1=\frac{\bigl(abcu,qabc/u;q\bigr)_{\infty}}
{\bigl(q,ab,ac,bc,au,qa/u,bu,qb/u,cu,qc/u;q\bigr)_{\infty}}
{}_3\phi_2\left(\begin{matrix} ab,ac,bc\\
abcu,qabc/u\end{matrix}\,; q,q\right).
\end{equation}
On the other hand, a straightforward residue computation shows that
\begin{equation*}
\begin{split}
L_2=-&\frac{\bigl(u^2/q^2;q\bigr)_{\infty}}
{\bigl(q,q,qa/u,qb/u,qc/u,au/q,bu/q,cu/q,u^2/q;q\bigr)_{\infty}}\\
&\qquad\times {}_7\phi_7\left(\begin{matrix}
q^2/u^2,q^2/u,-q^2/u,q,qa/u,qb/u,qc/u\\
q/u,-q/u,q^2/u^2,q^2/au,q^2/bu,q^2/cu,0\end{matrix}\,;
q,\frac{q^2}{abcu}\right).
\end{split}
\end{equation*}
Applying \cite[(3.8.9)]{GR} with the parameters $a,c,d,e$ and $f$
in \cite[(3.8.9)]{GR} replaced by $q^2/u^2$, $q$, $qa/u$, $qb/u$
and $qc/u$, respectively, we arrive at
\begin{equation*}
\begin{split}
L_2=-&\frac{\bigl(q/bc;q\bigr)_{\infty}\theta(u^2/q^2)}
{\bigl(q,q,qa/u,qb/u,qc/u,au/q,u^2/q;q\bigr)_{\infty}\theta(bu/q,cu/q)}\\
&\qquad\qquad\qquad\qquad\qquad\qquad
\times{}_3\phi_2\left(\begin{matrix} q/au,qb/u,qc/u\\
q^2/u^2,q^2/au\end{matrix}\,; q,\frac{q}{bc}\right).
\end{split}
\end{equation*}
Now applying the transformation formula \cite[(3.2.7)]{GR}
for ${}_3\phi_2$'s with the parameters $a,b,c,d$ and $e$ in
\cite[(3.2.7)]{GR} specialized to $q/au$, $qb/u$, $qc/u$,
$q^2/u^2$ and $q^2/au$, respectively, we obtain
\begin{equation}\label{stepB2}
\begin{split}
L_2=&-\frac{\bigl(q^2/abcu;q\bigr)_{\infty}\theta(u^2/q^2)}
{\bigl(q,qa/u,qb/u,qc/u,u^2/q;q\bigr)_{\infty}\theta(au/q,bu/q,cu/q)}\\
&\qquad\qquad\qquad\qquad\qquad\qquad
\times{}_3\phi_2\left(\begin{matrix} q/au,q/bu,q/cu\\
q^2/u^2,q^2/abcu\end{matrix}\,; q,q\right).
\end{split}
\end{equation}
Now combine \eqref{stepB1} and \eqref{stepB2}, and simplify
the terms using
\begin{equation}\label{theta}
\theta(qx^{-1})=\theta(x),\qquad \theta(qx)=(-x)^{-1}\theta(x),
\end{equation}
then we obtain
\begin{equation*}
\begin{split}
L_1+L_2=&\frac{1}{\bigl(q,qa/u,qb/u,qc/u,au,bu,cu;q\bigr)_{\infty}}\\
&\qquad\times\left\{\frac{\bigl(abcu,qabc/u;q\bigr)_{\infty}}
{\bigl(ab,ac,bc;q\bigr)_{\infty}}
{}_3\phi_2\left(\begin{matrix} ab,ac,bc\\
abcu, qabc/u\end{matrix}\,; q,q\right)\right.\\
&\quad\qquad\qquad\left.-\frac{q}{abcu}
\frac{\bigl(q^2/abcu,q^2/u^2;q\bigr)_{\infty}}
{\bigl(q/au,q/bu,q/cu;q\bigr)_{\infty}}
{}_3\phi_2\left(\begin{matrix}q/au,q/bu,q/cu\\
q^2/u^2,q^2/abcu\end{matrix}\,; q,q\right)\right\}.
\end{split}
\end{equation*}
By the three term transformation formula \cite[(3.3.1)]{GR}
for the ${}_3\phi_2$ basic hypergeometric series with the parameters 
$a,b,c,d$ and $e$ in \cite[(3.3.1)]{GR} specialized to
$qc/u$, $ac$, $bc$, $qabc/u$ and $qc/u$, respectively,
we obtain
\begin{equation*}
L_1+L_2=\frac{\bigl(qabc/u;q\bigr)_{\infty}\theta(abcu)}
{\bigl(q,ab,ac,bc,cu,qa/u,qb/u;q\bigr)_{\infty}\theta(au,bu)}
{}_2\phi_1\left(\begin{matrix} ac,bc\\
qabc/u\end{matrix}\,; q,\frac{q}{cu}\right).
\end{equation*}
Application of the $q$-Gauss sum \cite[(1.5.1)]{GR}
yields
\[ L_1+L_2=
\frac{1}{\bigl(q,ab,ac,bc;q\bigr)_{\infty}}
\frac{\theta(abcu)}{\theta(au,bu,cu)}.
\]
The evaluation of the one variable $q$-Macdonald-Mehta integral
\eqref{MacdonaldMehta} follows from this last formula by analytic
continuation.

\subsection{Second direct proof of Theorem \ref{MMint}}

For generic values of the parameters $a,b,c,u\in\mathbb{C}^\times$
satisfying $|a|,|b|,|c|<1$ and $q<|u|<1$, we can rewrite
the left hand side of \eqref{MacdonaldMehta} as
\begin{equation}\label{dversion}
\begin{split}
\frac{1}{2}\sum_{\stackrel{\scriptstyle{e\in \{a,b,c\}}}{k\in \mathbb{Z}_+}}
&\underset{x=eq^k}{\hbox{Res}}\left(\frac{\bigl(x^2,1/x^{2};q\bigr)_{\infty}}
{\bigl(ax,a/x,bx,b/x,cx,c/x;q\bigr)_{\infty}\theta(ux,u/x)x}\right)\\
+&\frac{1}{2}\sum_{k\in \mathbb{Z}}
\underset{x=uq^k}{\hbox{Res}}\left(\frac{\bigl(x^2,1/x^{2};q\bigr)_{\infty}}
{\bigl(ax,a/x,bx,b/x,cx,c/x;q\bigr)_{\infty}\theta(ux,u/x)x}\right)
\end{split}
\end{equation}
by shrinking the radius of the integration circle to zero while
picking up residues, cf. \cite[Section 4.10]{GR}.
The first three sums over $k\in \mathbb{Z}_+$ with fixed
$e\in \{a,b,c\}$ in \eqref{dversion} can be evaluated by the limit case
$d\rightarrow\infty$ in
Rogers' \cite[(2.7.1)]{GR} summation formula of a
very-well-poised ${}_6\phi_5$ series. For instance, the case $e=a$ yields
\begin{equation*}
\begin{split}
\sum_{k\in \mathbb{Z}_+}
&\underset{x=aq^k}{\hbox{Res}}\left(\frac{\bigl(x^2,1/x^{2};q\bigr)_{\infty}}
{\bigl(ax,a/x,bx,b/x,cx,c/x;q\bigr)_{\infty}\theta(ux,u/x)x}\right)\\
&\quad=\frac{\bigl(1/a^2;q\bigr)_{\infty}}
{\bigl(q,ab,b/a,ac,c/a;q\bigr)_{\infty}\theta(ua,u/a)}\,
{}_5\phi_5\left(\begin{matrix} a^2,qa,-qa,ab,ac\\
a,-a,qa/b,qa/c,0\end{matrix}\,; q,\frac{q}{bc}\right)\\
&\qquad\qquad=\frac{\bigl(q/bc;q\bigr)_{\infty}}
{\bigl(q,ab,ac;q\bigr)_{\infty}}
\frac{\theta(1/a^2)}{\theta(b/a,c/a,au,u/a)}.
\end{split}
\end{equation*}
The sums for $e=b,c$ can be obtained by interchanging the role
of $a$ and $e$ in the above formula. The fourth sum in
\eqref{dversion} (over $k\in \mathbb{Z}$)
can be evaluated by the limit case $e\rightarrow\infty$
in Bailey's \cite[(5.3.1)]{GR} summation formula of a
very-well-poised ${}_6\psi_6$ series, yielding
\begin{equation*}
\begin{split}
&\sum_{k\in \mathbb{Z}}\,
\underset{x=uq^k}{\hbox{Res}}\left(\frac{\bigl(x^2,1/x^{2};q\bigr)_{\infty}}
{\bigl(ax,a/x,bx,b/x,cx,c/x;q\bigr)_{\infty}\theta(ux,u/x)x}\right)\\
&\,=\frac{(1-1/u^{2})}
{\bigl(q,q,au,bu,cu,a/u,b/u,c/u;q\bigr)_{\infty}}\,
{}_5\psi_6\left(\begin{matrix} q/u,-q/u,a/u,b/u,c/u\\
1/u,-1/u,q/au,q/bu,q/cu,0\end{matrix}\,; q,\frac{q}{abcu}\right)\\
&\qquad\qquad=\frac{\bigl(q/ab,q/ac,q/bc;q\bigr)_{\infty}}
{\bigl(q;q\bigr)_{\infty}}
\frac{\theta(1/u^2)}{\theta(au,bu,cu,a/u,b/u,c/u)}.
\end{split}
\end{equation*}
Set $\chi=\log(u)$, then it remains to evaluate
\begin{equation*}
\begin{split}
f(\chi)=&\frac{\bigl(q/bc;q\bigr)_{\infty}}
{\bigl(q,ab,ac;q\bigr)_{\infty}}
\frac{\theta(1/a^2)}{\theta(b/a,c/a,ae^{\chi},e^{\chi}/a)}\\
&+\frac{\bigl(q/ac;q\bigr)_{\infty}}
{\bigl(q,ab,bc;q\bigr)_{\infty}}
\frac{\theta(1/b^2)}{\theta(a/b,c/b,be^{\chi},e^{\chi}/b)}\\
&+\frac{\bigl(q/ab;q\bigr)_{\infty}}
{\bigl(q,bc,ac;q\bigr)_{\infty}}
\frac{\theta(1/c^2)}{\theta(b/c,a/c,ce^{\chi},e^{\chi}/c)}\\
&+\frac{\bigl(q/ab,q/ac,q/bc;q\bigr)_{\infty}}
{\bigl(q;q\bigr)_{\infty}}
\frac{\theta(e^{-2\chi})}{\theta(ae^{\chi},be^{\chi},
ce^{\chi},ae^{-\chi},be^{-\chi},ce^{-\chi})},
\end{split}
\end{equation*}
which we consider as a meromorphic function in $\chi\in
\mathbb{C}$, with fixed (generic) parameters $a,b,c$.
We consider first the meromorphic function
\[g(\chi)=\theta(ae^{\chi},be^{\chi},ce^{\chi})f(\chi).
\]
Observe that the possible poles of $g$ are located at $\log(e)+\Lambda$
for $e=a,b,c$, where $\Lambda\subset\mathbb{C}$ is the lattice
\[\Lambda=\mathbb{Z}\,\log(q)+\mathbb{Z}\,2\pi i.
\]
Furthermore, the poles are at most simple for
generic parameter values. We show now that $g$ is in fact analytic.
Observe that $g$ is $2\pi i$-periodic. Furthermore, by \eqref{theta},
\[g(\chi+\log(q))=\left(\frac{-1}{abce^{\chi}}\right)g(\chi),
\]
i.e. $g$ is quasi-periodic with quasi-period $\log(q)$.
So $g$ is quasi-periodic with respect to the period lattice $\Lambda$.
In view of the symmetry of $g$ in the parameters $a,b$ and $c$, we
conclude that $g$ is analytic if the residue of $g(\chi)$ at
$\alpha:=\log(a)$ is zero. This follows from the observation that
\begin{equation*}
\begin{split}
\lim_{\chi\rightarrow \alpha}(1-e^{\chi-\alpha})g(\chi)=
&\frac{\bigl(q/bc;q\bigr)_{\infty}}
{\bigl(q,q,q,ab,ac;q\bigr)_{\infty}}\frac{\theta(1/a^2,ab,ac)}
{\theta(b/a,c/a)}\\
&-\frac{\bigl(q/ab,q/ac,q/bc;q\bigr)_{\infty}}
{\bigl(q,q,q;q\bigr)_{\infty}}\frac{\theta(1/a^2)}
{\theta(b/a,c/a)}\\
&=0.
\end{split}
\end{equation*}
We conclude that the function
\begin{equation*}
\begin{split}
h(\chi)=\frac{g(\chi)}{\theta(abce^{\chi})}&=
\frac{\theta(ae^{\chi},bq^{\chi},ce^{\chi})}
{\theta(abce^{\chi})}f(\chi)\\
&=\frac{\bigl(q/bc;q\bigr)_{\infty}}
{\bigl(q,ab,ac;q\bigr)_{\infty}}
\frac{\theta(1/a^2,be^{\chi},ce^{\chi})}
{\theta(b/a,c/a,e^{\chi}/a,abce^{\chi})}\\
&+\frac{\bigl(q/ac;q\bigr)_{\infty}}
{\bigl(q,ab,bc;q\bigr)_{\infty}}
\frac{\theta(1/b^2,ae^{\chi},ce^{\chi})}
{\theta(a/b,c/b,e^{\chi}/b,abce^{\chi})}\\
&+\frac{\bigl(q/ab;q\bigr)_{\infty}}
{\bigl(q,bc,ac;q\bigr)_{\infty}}
\frac{\theta(1/c^2,ae^{\chi},be^{\chi})}{\theta(b/c,a/c,
e^{\chi}/c,abce^{\chi})}\\
&+\frac{\bigl(q/ab,q/ac,q/bc;q\bigr)_{\infty}}
{\bigl(q;q\bigr)_{\infty}}
\frac{\theta(e^{-2\chi})}{\theta(abce^{\chi},
ae^{-\chi},be^{-\chi},ce^{-\chi})}
\end{split}
\end{equation*}
defines an elliptic function with respect to the period lattice
$\Lambda$,
with at most one pole in a fundamental domain of $\mathbb{C}/\Lambda$.
Thus $h$ is constant, so in particular,
\[f(\chi)=\frac{\theta(abce^{\chi})}{\theta(ae^{\chi},be^{\chi},
ce^{\chi})}h(-\alpha), \qquad \alpha=\log(a).
\]
By the explicit expression for $h$ we have
\[h(-\alpha)=\frac{2}{\bigl(q,ab,ac,bc;q\bigr)_{\infty}},
\]
hence the left hand side of \eqref{MacdonaldMehta} is equal to
\[ \frac{1}{2}f(\log(u))=\frac{1}{\bigl(q,ab,ac,bc;q\bigr)_{\infty}}
\frac{\theta(abcu)}{\theta(au,bu,cu)},
\]
which completes the proof of Theorem \ref{MMint}.

%%%%%%%%%%%%%%%%%%%%%%%%%%%%%%%%%%%%%%%%%%%%%%%%%%%%%%%%%%%%%%%%
%%                                                            %%
%%                       Bibliography                         %%
%%                                                            %%
%%%%%%%%%%%%%%%%%%%%%%%%%%%%%%%%%%%%%%%%%%%%%%%%%%%%%%%%%%%%%%%%

\bibliographystyle{amsplain}

\end{document}